# Integrals and Series Related to the Surface Area of Arbitrary Ellipsoids

Richard A. Krajcik and Kelly D. McLenithan

Los Alamos National Laboratory

## Abstract

Integrals related to the surface area of arbitrary ellipsoids, $\frac{x^2}{a^2} + \frac{y^2}{b^2} + \frac{z^2}{c^2} = 1$, are derived, evaluated, and compared with each other and existing integrals found in the literature for $a > b > c > 0$ and $c > b > a > 0$. For ellipsoids where $a > b > c > 0$, six integrals are derived using standard integration techniques and evaluated using the known surface area expression for an arbitrary ellipsoid. The six integrals are not found in the literature in either ellipsoid-specific or general form. For $a > b > c > 0$, infinite series representations are derived for two ellipsoid-specific integrals. With one exception, when $c > b > a > 0$, integrals can be evaluated directly using integral tables found in the literature, including integrals that lead to a confirmation of the known expression of the surface area of an arbitrary ellipsoid with $a$ and $c$ interchanged. Analytic extensions of integrals evaluated for $a > b > c > 0$ are then compared with integrals evaluated for $c > b > a > 0$ and with integrals found in the literature. While some integrals can be found in the literature, three additional integrals do not appear to have a counterpart in the literature. Nine integrals derived and evaluated in this study, as expressed in their general form, are expected to have application beyond that of the surface area of ellipsoids. A brief list of historical references is provided and some recent communications regarding the surface area of ellipsoids are discussed in the context of this study.

## 1    Introduction

An exact, closed-form expression for the surface area of an arbitrary ellipsoid has been known (at least) since the time of Legendre (1752-1833). Legendre's expression for the surface area of an arbitrary ellipsoid is presented below from *Exercices de calcul intégral* (p. 191), published (Huzard-Courier, Paris) in 1811, namely,

(1.1) $\quad S(a,b,c) = 2\pi c^2 + \frac{2\pi ab}{\sin v}\left[ \frac{c^2}{a^2} F(v,b') + \frac{a^2 - c^2}{a^2} E(v,b') \right].$     **LE1 191**

Here $S(a,b,c)$ represents the total surface area of an ellipsoid with ordered variables $(a,b,c)$ such that $a > b > c > 0$, $\cos v = \frac{c}{a}$, $b'^2 = \frac{b^2 - c^2}{b^2 \sin^2 v}$, and $F(v,b')$ and $E(v,b')$ are incomplete elliptic integrals of the first and second kind, respectively, with argument $v$ and modulus $b'$. Note that we have reversed the order of argument and modulus for the incomplete elliptic integrals in Eq.

(1.1). In the same publication, Legendre gave expressions for the surface area of an oblate spheroid,

(1.2) $\quad S(a,a,c) = 2\pi a^2 + \dfrac{\pi c^2}{\sin \nu} \ln\left(\dfrac{1+\sin \nu}{1-\sin \nu}\right),$ **LE1 191**

where $a = b > c > 0$, and for a prolate spheroid,

(1.3) $\quad S(a,b,b) = \dfrac{2\pi ab}{\sin \nu}[\nu + \sin \nu \cos \nu],$ **LE1 191**

where $a > b = c > 0$. It is remarkable that Legendre published Eqs. (1.1), (1.2), and (1.3) nearly two centuries ago.

In his 1811 publication (above), Legendre refers to an earlier work by Monge (1746-1818), *Feuilles d'analyse appliquée à la géometrie,* first presented as lectures at the École Polytechnique (Paris) in 1795, later published in 1801, and still later expanded in *Application de l'analyse à la géometrie* published in 1807. These (and other) publications by Monge provided a foundation for descriptive and analytic geometry. While the work of Monge is fundamental to an understanding of descriptive and analytic geometry, the expression for the surface area of an arbitrary ellipsoid given by Legendre in 1811 appears to be the earliest published expression for the surface area of an ellipsoid in terms of incomplete elliptic integrals that can be easily recognized today. Unfortunately, since Legendre first published his expression for the surface area of an ellipsoid in terms of elliptic integrals in 1811, the expression seems to have been alternately forgotten, rediscovered, forgotten again, rediscovered again, and subsequently rederived a number of times during the interim period of nearly two centuries. An abbreviated chronology of the surface area of an arbitrary ellipsoid is provided in the general references[1-22], together with some additional comments intended to clarify the literature. It is important to note that the expression for the surface area of an arbitrary ellipsoid has not always been published, nor clearly represented when published, even in very recent standard references[16, 21]. Of the many references that provide an exact, closed-form expression for the surface area of an arbitrary ellipsoid, we shall most closely follow the work of Bowman[7, 9].

In 1953 Bowman[7] published an elegant derivation of the surface area of an ellipsoid in a book entitled *Introduction to Elliptic Functions with Applications.* Bowman's 1953 book was subsequently republished in 1961[9], but both books have since gone out of print. In his derivation, Bowman assumed that the ellipsoid was characterized by

(1.4) $\quad \dfrac{x^2}{a^2} + \dfrac{y^2}{b^2} + \dfrac{z^2}{c^2} = 1,$ **BO 31 (30)**

with eccentricities given by

(1.5) $\quad e_1^2 = \dfrac{a^2 - c^2}{a^2}$ and $e_2^2 = \dfrac{b^2 - c^2}{b^2},$ **BO 32**

where $a^2 > b^2 > c^2$. Here $a$, $b$, and $c$ are assumed to be positive, real numbers, with positive square roots for the eccentricities $e_1$ and $e_2$. Since $a > b > c > 0$, it follows that $1 > e_1 > e_2 > 0$. Note that we have adopted the notation of Gradshteyn and Ryzhik[23] for functions wherever possible, and also for bibliographic references, so that **BO 31 (30)** refers to Bowman (1961), p. 31, Eq. (30). In what follows, we shall show that the relative magnitude of the parameters $a$, $b$, and $c$ is critical in the functional characterization of integrals related to the surface area of an arbitrary ellipsoid.

As originally published, Bowman's expression for the surface area of an arbitrary ellipsoid was given in terms of elliptic functions, namely,

$$(1.6) \quad S = 2\pi c^2 + \frac{2\pi b}{\sqrt{a^2 - c^2}}\left[(a^2 - c^2)E(\theta) + c^2\theta\right], \qquad \textbf{BO 32 (40)}$$

where

$$(1.7) \quad k^2 = \frac{e_2^2}{e_1^2} \text{ and } e_1 = \text{sn}(\theta, k). \qquad \textbf{BO 32 (37, 38)}$$

In terms of elliptic integrals, Eq. (1.6) can be rewritten as

$$(1.8) \quad S(a,b,c) = 2\pi c^2 + \frac{2\pi b}{\sqrt{a^2 - c^2}}\left[(a^2 - c^2)E(\varphi, k) + c^2 F(\varphi, k)\right],$$

where $a > b > c > 0$, and

$$(1.9) \quad F(\varphi, k) = \theta \text{ and } E(\varphi, k) = E(\theta) \qquad \textbf{BO 16 (8, 9)}$$

are incomplete elliptic integrals of the first and second kind, respectively, with argument $\varphi$ and modulus $k$. With $\text{sn}\,\theta \equiv \text{sn}(\theta, k) = \sin\varphi$ [**BO 9, 16**], we have from Eqs. (1.5) and (1.7) above,

$$(1.10) \quad \varphi = \arcsin e_1 = \arcsin\sqrt{\frac{a^2 - c^2}{a^2}} \text{ and } k = \frac{e_2}{e_1} = \sqrt{\frac{a^2}{b^2}\left(\frac{b^2 - c^2}{a^2 - c^2}\right)}.$$

Note that we have used the notation of Gradshteyn and Ryzhik[23] for incomplete elliptic integrals in Eqs. (1.8) and (1.9) instead of the notation originally provided by Bowman, namely $F(k, \varphi) = \theta$ and $E(k, \varphi) = E(\theta)$. Note too that Bowman's function $E(\theta)$ represents an incomplete elliptic integral of the second kind, $E(\varphi, k)$, expressed in elliptic function notation, and not a complete elliptic integral of the second kind with modulus $\theta$ [See for example, **BO 16 (9), BY 8 (110.03), and ER2 343 (12)**]. Finally, note that we shall use the notation of Gradshteyn and Ryzhik[23] for complete elliptic integrals, so that complete elliptic integrals of the first and second kind, $\boldsymbol{F}(k)$ and $\boldsymbol{E}(k)$, respectively, are distinguished from other elliptic integral forms by their bold type and single variable. The complementary modulus $k' = \sqrt{1-k^2}$ is defined in the usual way [**BO 9 (7)**]. Since $1 > e_1 > e_2 > 0$, we have $\frac{\pi}{2} > \varphi > 0$ (principal value), and $1 > k > 0$ and $1 > k' > 0$.

If $a$ and $b$ are expressed in terms of $c$ and the eccentricities, then the surface area of the ellipsoid can be written as

$$(1.11) \quad S(a,b,c) = 2\pi c^2 \left\{ 1 + \sqrt{\frac{1-e_1^2}{1-e_2^2}} \left[ \frac{e_1^2 E(\varphi,k) + (1-e_1^2)F(\varphi,k)}{e_1(1-e_1^2)} \right] \right\},$$

where

$$(1.12) \quad a = \frac{c}{\sqrt{1-e_1^2}} \text{ and } b = \frac{c}{\sqrt{1-e_2^2}},$$

while $\varphi$ and $k$ are given by Eq. (1.10). Note that $e_1$ and $e_2$ contain the ellipsoid shape information, while the lengths of $a$ and $b$, and the surface area $S(a,b,c)$ are scaled by $c$ and $c^2$, respectively. An illustration of an ellipsoid where $a > b > c > 0$ is presented in Fig. 1. Here $(a,b,c)$ have been arbitrarily set to $(2.0, 1.5, 1.0)$.

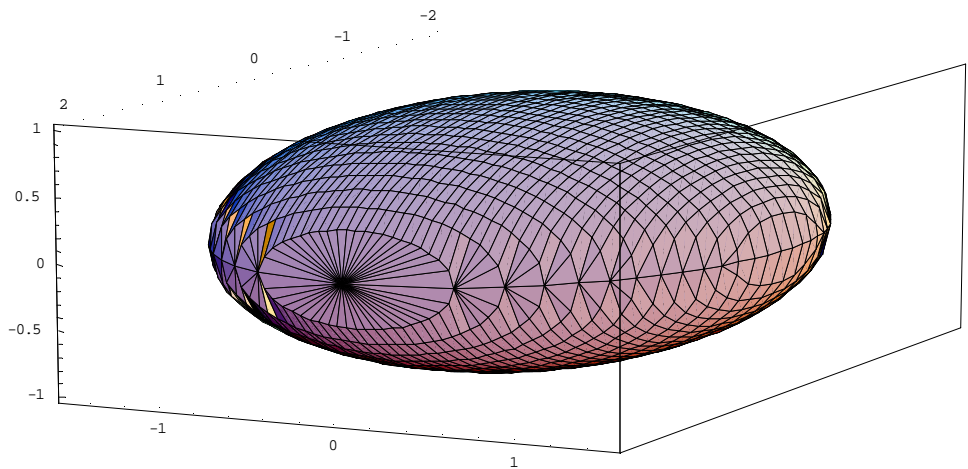

Figure 1. An Illustration of an Ellipsoid where $a(2.0) > b(1.5) > c(1.0) > 0$.

If the surface area of the ellipsoid illustrated in Fig. 1 is projected onto the x-y plane, as indicated in Fig. 2, then it would seem to be a straightforward matter to integrate the simple equation that characterizes an ellipsoid where $(a,b,c)$ have been chosen such that $a > b > c > 0$. However, as the subject of this study already implies, integration is not nearly as straightforward as it might be for an object as simple as an ellipsoid where $a > b > c > 0$. Without resorting to elliptic functions [**BO 32 (39)**] or complex variables, we shall show in Section 2 that none of the integrals that arise from a straightforward application of integration techniques lead directly to the known surface area expression as given by either Eq. (1.8) or (1.11). However, there is a positive side to this difficulty, namely that all the integrals considered in Section 2 are "one integral short" from reproducing the known surface area expression for an ellipsoid; so that the surface area expression itself can be used to evaluate the integrals in the case of an arbitrary ellipsoid where $a > b > c > 0$. The general form of these integrals, expected to have application beyond the surface area of an arbitrary ellipsoid, could not be found in the literature.

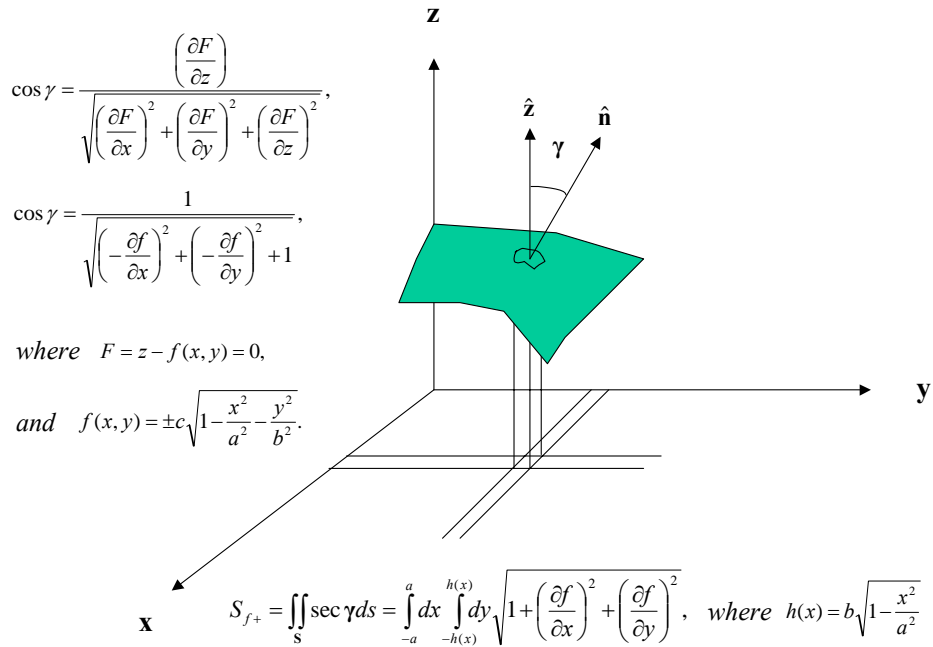

Figure 2. Projection of Ellipsoid Surface Area onto the x-y Plane.

Next consider an ellipsoid where $a > b > c > 0$ (Fig. 1) rotated 90° about the y-axis. This rotation has the effect of reversing the order of the relative magnitude of $(a,b,c)$, so that in the rotated case, $c > b > a > 0$. Figure 3 illustrates a rotated ellipsoid where $c > b > a > 0$, and $(a,b,c)$ have been set to $(1.0, 1.5, 2.0)$. If the surface area of the rotated ellipsoid is again projected onto the x-y plane (Fig. 2), then this relatively minor difference in orientation allows some integrals to be evaluated such that they directly reproduce the rotated surface area result as given by Eq. (1.13) below. Here interchanging $a$ and $c$ in Eq. (1.8) produces the desired result for the rotated ellipsoid,

(1.13) $\quad S(c,b,a) = 2\pi a^2 + \dfrac{2\pi b}{\sqrt{c^2 - a^2}} \left[ (c^2 - a^2) E(\overline{\varphi}, \overline{k}) + a^2 F(\overline{\varphi}, \overline{k}) \right],$

(1.14)  $\bar{\varphi} = \arcsin\sqrt{\dfrac{c^2 - a^2}{c^2}}$ and $\bar{k} = \sqrt{\dfrac{c^2}{b^2}\left(\dfrac{b^2 - a^2}{c^2 - a^2}\right)}$,

and where $\bar{\varphi}$ and $\bar{k}$ denote the argument and modulus of the elliptic integrals $F(\bar{\varphi},\bar{k})$ and $E(\bar{\varphi},\bar{k})$ for the rotated ellipsoid where $c > b > a > 0$. If the eccentricities of the rotated ellipsoid are denoted by $\bar{e}_1$ and $\bar{e}_2$, and defined with respect to the semi-major axis, $c$, then

(1.15)  $\bar{e}_1^2 = \dfrac{c^2 - a^2}{c^2}$ and $\bar{e}_2^2 = \dfrac{c^2 - b^2}{c^2}$, where $c^2 > b^2 > a^2$.

Note that the definition of $\bar{e}_2$ is not what one might expect from the interchange of $a$ and $c$ in Eq. (1.5). However, the definitions of eccentricities given by Eq. (1.15), and subsequently $(\bar{f}_1, \bar{f}_2)$ given by Eq. (1.18), allow the surface area Eq. (3.3) to separate for $c > b > a > 0$ in the same way that the eccentricities defined by Eq. (1.5) allow the surface area Eq. (2.3) to separate when $a > b > c > 0$.

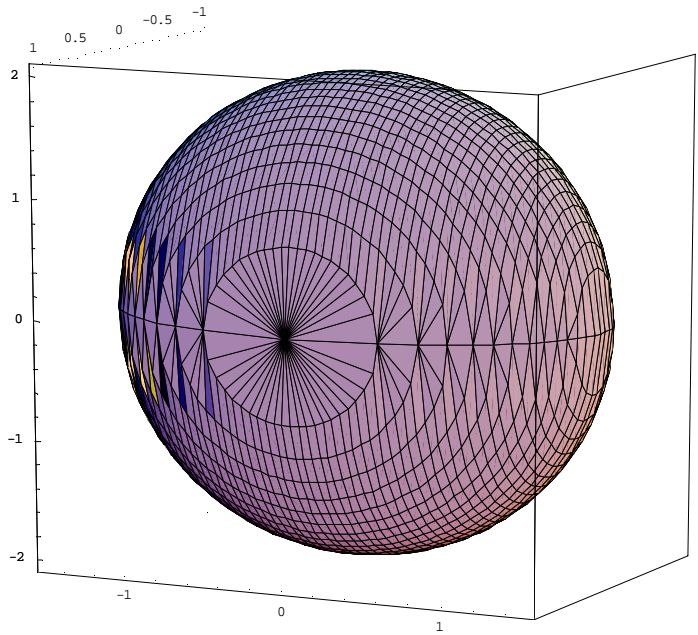

Figure 3. An Illustration of an Ellipsoid where $c(2.0)>b(1.5)>a(1.0)>0$.

With $\bar{e}_1$ and $\bar{e}_2$ so defined, Eqs. (1.13) and (1.14) can be rewritten in terms of $a$ and the eccentricities as,

(1.16)  $S(c,b,a) = 2\pi a^2 \left\{ 1 + \sqrt{\dfrac{1-\bar{e}_2^2}{1-\bar{e}_1^2}} \left[ \dfrac{\sqrt{1-\bar{e}_1^2}\,F(\bar{\varphi},\bar{k})}{\bar{e}_1} + \dfrac{\bar{e}_1 E(\bar{\varphi},\bar{k})}{\sqrt{1-\bar{e}_1^2}} \right] \right\}$,

where

(1.17) $\bar{\varphi} = \arcsin \bar{e}_1$ and $\bar{k} = \sqrt{\dfrac{\bar{e}_1^2 - \bar{e}_2^2}{\bar{e}_1^2(1-\bar{e}_2^2)}} = \sqrt{1 - \dfrac{\left(\dfrac{\bar{e}_2^2}{1-\bar{e}_2^2}\right)}{\left(\dfrac{\bar{e}_1^2}{1-\bar{e}_1^2}\right)}}$ .

Again, we assume the positive square root for the eccentricities. Since $c > b > a > 0$, it follows that $1 > \bar{e}_1 > \bar{e}_2 > 0$. Now define $\bar{f}_1$ and $\bar{f}_2$ such that

(1.18) $\bar{f}_1 = \dfrac{\bar{e}_1}{\sqrt{1-\bar{e}_1^2}} = \sqrt{\dfrac{c^2 - a^2}{a^2}}$ and $\bar{f}_2 = \dfrac{\bar{e}_2}{\sqrt{1-\bar{e}_2^2}} = \sqrt{\dfrac{c^2 - b^2}{b^2}}$ .

With $\bar{f}_1$ and $\bar{f}_2$ so defined, Eqs. (1.16) and (1.17) can be rewritten as

(1.19) $S(c,b,a) = 2\pi a^2 \left\{ 1 + \sqrt{\dfrac{1+\bar{f}_1^2}{1+\bar{f}_2^2}} \left[ \dfrac{F(\bar{\varphi},\bar{k})}{\bar{f}_1} + \bar{f}_1 E(\bar{\varphi},\bar{k}) \right] \right\}$ ,

where

(1.20) $\bar{\varphi} = \arctan \bar{f}_1$ and $\bar{k} = \sqrt{1 - \dfrac{\bar{f}_2^2}{\bar{f}_1^2}}$ ,                     **GR 55 (1.624.7)**

and

(1.21) $a = \dfrac{c}{\sqrt{1+\bar{f}_1^2}}$ and $b = \dfrac{c}{\sqrt{1+\bar{f}_2^2}}$ .

Since $1 > \bar{e}_1 > \bar{e}_2 > 0$, it follows that $\infty > \bar{f}_1 > \bar{f}_2 > 0$, and it follows that $\dfrac{\pi}{2} > \bar{\varphi} > 0$ and $1 > \bar{k} > 0$.

In Section 3, we show that $\bar{f}_1$ and $\bar{f}_2$ provide a natural description of surface area integrals for an ellipsoid where $c > b > a > 0$. Note that the difference between surface area integrals for an ellipsoid, both rotated ($c > b > a > 0$) and unrotated ($a > b > c > 0$), is similar to the difference found between surface area integrals for prolate and oblate spheroids. In short, the shape of the spheroid determines the condition of integration, which in turn determines the functions that characterize the integral.

Note that the surface area of an ellipsoid, $a > b > c > 0$, given by Eq. (1.11), has the correct limit when $a = b = r > c > 0$, namely,

$$(1.22) \quad S(r,r,c) \equiv \lim_{a \to r} S(a,r,c) = 2\pi c^2 \left( 1 + \frac{e_1 E(\varphi,1)}{1-e_1^2} + \frac{F(\varphi,1)}{e_1} \right),$$

where

$$(1.23) \quad \begin{Bmatrix} F(\varphi,1) \\ E(\varphi,1) \end{Bmatrix} = \begin{Bmatrix} \ln(\tan\varphi + \sec\varphi) \\ \sin\varphi \end{Bmatrix}, \qquad \text{BY 10 (111.04)}$$

and

$$(1.24) \quad \varphi = \arcsin e_1 \text{ and } e_1 = \sqrt{\frac{r^2 - c^2}{r^2}} = e_2.$$

This allows the expression for the surface area of an oblate spheroid to be recovered, namely,

$$(1.25) \quad S(r,r,c) = 2\pi r^2 + \frac{\pi r c^2}{\sqrt{r^2 - c^2}} \ln\left( \frac{r + \sqrt{r^2 - c^2}}{r - \sqrt{r^2 - c^2}} \right). \qquad \text{ZW 364 (4.18.1) (r>c)}$$

Here the oblate spheroid is symmetric in the x-y plane with radius $r$ and $r > c > 0$.

Similarly, the surface area of an ellipsoid, $c > b > a > 0$, given by Eq. (1.19), has the correct limit when $c > b = a = r > 0$, namely,

$$(1.26) \quad S(c,r,r) \equiv \lim_{a \to r} S(c,r,a) = 2\pi r^2 \left( 1 + \frac{F(\bar{\varphi},0)}{\bar{f}_1} + \bar{f}_1 E(\bar{\varphi},0) \right),$$

where

$$(1.27) \quad \begin{Bmatrix} F(\bar{\varphi},0) \\ E(\bar{\varphi},0) \end{Bmatrix} = \begin{Bmatrix} \bar{\varphi} \\ \bar{\varphi} \end{Bmatrix}, \qquad \text{BY 10 (111.01)}$$

and

$$(1.28) \quad \bar{\varphi} = \arctan \bar{f}_1 \text{ and } \bar{f}_1 = \sqrt{\frac{c^2 - r^2}{r^2}} = \bar{f}_2,$$

so that the expression for the surface area of a prolate spheroid is recovered, namely,

$$(1.29) \quad S(c,r,r) = 2\pi r^2 + \frac{2\pi r c^2}{\sqrt{c^2 - r^2}} \arcsin \sqrt{\frac{c^2 - r^2}{c^2}}, \qquad \text{ZW 364 (4.18.1) (c>r)}$$

Here the prolate spheroid is symmetric in the x-y plane with radius $r$ and $c > r > 0$.

With the ellipsoid oriented relative to the x-y plane such that $c > b > a > 0$, integrals could be found in the literature that directly lead to a confirmation of the known expression for the surface area. However, while most integrals related to the surface area of an ellipsoid where $c > b > a > 0$ could be evaluated directly using integral tables found in the literature, one integral could not be directly evaluated. In that case, the same approach described in Section 2 was employed to evaluate a surface area integral for an ellipsoid where $c > b > a > 0$; that is, Eq. (1.13) or (1.19) was used to evaluate the derived integral. Integrals related to the surface area of an ellipsoid where $c > b > a > 0$ are discussed and tabulated by integral identifier $I_M(c,b,a)$ in Section 3.

In Section 4, the analytic extension of integrals related to the surface area of an ellipsoid, derived and evaluated in Sections 2 and 3, are discussed and compared with each other and with integrals found in the literature. When we extend the range of the argument $\varphi$ and modulus $k$ of elliptic integrals, we shall use the expressions and adopt the notation of Byrd and Friedman[24, 25]. Naturally, surface area integrals evaluated in Section 2 for an ellipsoid where $a > b > c > 0$ can be expressed as analytic extensions of the corresponding surface area integrals evaluated in Section 3 for an ellipsoid where $c > b > a > 0$. This is similar to the situation with prolate and oblate spheroids. In addition, we note that some surface area integrals evaluated in Sections 2 and 3 are analytic extensions of integrals found in the literature. However, of the twelve integrals and extensions related to the surface area of an ellipsoid, derived and evaluated in Sections 2 and 3 and extended in Section 4, only three could be found in the literature. The remaining nine integrals are perhaps the most interesting contribution from this study of integrals and series related to the surface area of arbitrary ellipsoids.

Finally, we discuss some references and communications related to the surface area of an ellipsoid in the literature, including a communication by Dieckmann[19], where Dieckmann applied Mathematica 4 and used two integrals from Gradshteyn and Ryzhik[23]; namely **GR 625 (6.113.2)** and **GR 626 (6.123)**, to derive the surface area of an ellipsoid.

## 2    Ellipsoids (a,b,c) where a>b>c>0

The surface area of an ellipsoid $(a > b > c > 0)$ projected onto the x-y plane can be written as

$$(2.1) \quad S(a,b,c) = 8 \int_0^a \int_0^{h(x)} \sqrt{1 + \left(\frac{\partial f}{\partial x}\right)^2 + \left(\frac{\partial f}{\partial y}\right)^2}\, dy\, dx,$$

where

$$(2.2) \quad h(x) = b\sqrt{1 - \frac{x^2}{a^2}} \quad \text{and} \quad f(x,y) = c\sqrt{1 - \frac{x^2}{a^2} - \frac{y^2}{b^2}},$$

and the eight-fold symmetry of the ellipsoid has been used to reduce the range of integration. Figures 1 and 2 illustrate an ellipsoid where $a > b > c > 0$ with its surface area projected onto the x-y plane. Using Eqs. (1.5) and (2.2), Eq. (2.1) can be rewritten as

$$(2.3) \quad S(a,b,c) = 8\int_0^a \int_0^{h(x)} \frac{\sqrt{1 - \frac{e_1^2 x^2}{a^2} - \frac{e_2^2 y^2}{b^2}}}{\sqrt{1 - \frac{x^2}{a^2} - \frac{y^2}{b^2}}} \, dy\, dx,$$

then simplified with scaled elliptic variables where

$$(2.4) \quad x = a\hat{x} \text{ and } y = b\hat{y},$$

so that,

$$(2.5) \quad S(a,b,c) = 8ab \int_0^1 \int_0^{\sqrt{1-\hat{x}^2}} \frac{\sqrt{1 - e_1^2 \hat{x}^2 - e_2^2 \hat{y}^2}}{\sqrt{1 - \hat{x}^2 - \hat{y}^2}} \, d\hat{y}\, d\hat{x}.$$

Further simplification is possible with scaled polar coordinates, where

$$(2.6) \quad \hat{x} = \rho \cos\phi \text{ and } \hat{y} = \rho \sin\phi,$$

so that,

$$(2.7) \quad S(a,b,c) = 8ab \int_0^1 \int_0^{\pi/2} \frac{\sqrt{1 - (e_1^2 \cos^2\phi + e_2^2 \sin^2\phi)\rho^2}}{\sqrt{1 - \rho^2}} \rho\, d\phi\, d\rho,$$

where $\frac{\pi}{2} \geq \phi \geq 0$ and $1 \geq \rho \geq 0$. Scaled elliptic variables reduce the x-y ellipse to a pseudo-circle, while scaled orthogonal polar coordinates,

$$(2.8) \quad g_{\rho\phi} = \left(\frac{\partial \hat{x}}{\partial \rho}\right)\left(\frac{\partial \hat{x}}{\partial \phi}\right) + \left(\frac{\partial \hat{y}}{\partial \rho}\right)\left(\frac{\partial \hat{y}}{\partial \phi}\right) = g_{\phi\rho} = 0, \qquad \textbf{GR 1042 (10.51)}$$

simplify the differential area such that $d\hat{y}\, d\hat{x} = \rho\, d\phi\, d\rho$. Equation (2.7) could also be obtained by the use of non-orthogonal, generalized polar coordinates [**PR1 565 (3.1.1.2)**], but only because the differential surface area contains no cross terms (for non-orthogonal, generalized polar coordinates, $g_{\rho\phi} \neq 0$).

With a change of variables, namely $\rho = \sin\theta$, Eq. (2.7) can be written in a familiar form,

$$(2.9) \quad S(a,b,c) = 8ab \int_0^{\pi/2} \int_0^{\pi/2} \sin\theta \sqrt{1 - q^2 \sin^2\theta} \, d\theta\, d\phi,$$

where

(2.10) $\quad q^2 = e_1^2 \cos^2\phi + e_2^2 \sin^2\phi$.

Note that $1 > e_1 > e_2 > 0$ and $\frac{\pi}{2} \geq \phi \geq 0$ imply that $1 > q^2 > 0$. Equation (2.9) can be equally written as

(2.11) $\quad S(a,b,c) = 8ab \int_0^{\pi/2} \int_0^{\pi/2} \sin\theta \sqrt{1-e_1^2 \sin^2\theta} \sqrt{1+w^2 \sin^2\phi} \, d\theta d\phi$,

where

(2.12) $\quad w^2 = \dfrac{(e_1^2 - e_2^2)\sin^2\theta}{1 - e_1^2 \sin^2\theta}$.

Note that $1 > e_1 > e_2 > 0$ and $\frac{\pi}{2} \geq \theta \geq 0$ imply that $\infty > w^2 > 0$. Because the limits of $\theta$ and $\phi$ do not depend upon one another, Eqs. (2.9) and (2.11) can be integrated in either order. However, before doing so, it is instructive to rewrite Eq. (2.9) in terms of $(a,b,c)$. Applying Eq. (1.5) to Eq. (2.9) results in

(2.13) $\quad S(a,b,c) = 8 \int_0^{\pi/2} \int_0^{\pi/2} \sin\theta \sqrt{b^2 c^2 \sin^2\theta \cos^2\phi + a^2 c^2 \sin^2\theta \sin^2\phi + a^2 b^2 \cos^2\theta} \, d\theta d\phi$,

which is the same ellipsoid surface area integral as that provided by Mathematica 4 and 5[26] in the "Formula Gallery"[17] and that provided in a communication by Cantrell[18]. Furthermore, Eq. (2.13) provides the starting point for the Mathematica 4 and Gradshteyn and Ryzhik[23] derivation of the surface area of an ellipsoid by Dieckmann[19].

### 2.1 The Integral $I_1(a,b,c) = \int_0^\alpha \dfrac{u E(u) du}{\left(k'^2 + k^2 u^2\right)^2 \sqrt{\alpha^2 - u^2}}$

Begin the process by integrating Eq. (2.11), starting with $\phi$, so that,

(2.14) $\quad S(a,b,c) = 8ab \int_0^{\pi/2} \sin\theta \sqrt{1-e_1^2 \sin^2\theta} \left[ \sqrt{1+w^2} E(\alpha, \dfrac{w}{\sqrt{1+w^2}}) \right]_0^{\pi/2} d\theta$, **GR 199 (2.597.2)**

where $\alpha = \arcsin\left(\dfrac{\sqrt{1+w^2}\sin\phi}{\sqrt{1+w^2\sin^2\phi}}\right)$ and $w^2$ is given by Eq. (2.12). With these definitions, Eq. (2.14) can be reduced to

$$(2.15) \quad S(a,b,c) = 8ab \int_0^{\pi/2} \sin\theta\sqrt{1-e_2^2\sin^2\theta}\, E\left(\dfrac{\sqrt{e_1^2-e_2^2}\sin\theta}{\sqrt{1-e_2^2\sin^2\theta}}\right)d\theta,$$

where $E(k) \equiv E\left(\dfrac{\pi}{2}, k\right)$ is a complete elliptic integral of the second kind [**GR 852 (8.112.2)**]. To rewrite Eq. (2.15) in a more conventional form, let $u = \dfrac{\sqrt{e_1^2-e_2^2}\sin\theta}{\sqrt{1-e_2^2\sin^2\theta}}$, then, after a little algebra,

$$(2.16) \quad S(a,b,c) = 8ab \int_0^{\sqrt{\frac{e_1^2-e_2^2}{1-e_2^2}}} \dfrac{(e_1^2-e_2^2)^{\frac{3}{2}} u E(u) du}{\left(e_1^2-e_2^2+e_2^2 u^2\right)^2 \sqrt{(e_1^2-e_2^2)-(1-e_2^2)u^2}}.$$

Since $E(u) = \dfrac{\pi}{2} F(-\dfrac{1}{2},\dfrac{1}{2};1;u^2)$ [**GR 853 (8.114.1)**], where $F(\alpha,\beta;\gamma,z)$ is the Gauss hypergeometric series, and $F(-\dfrac{1}{2},\dfrac{1}{2};1;u^2)$ converges throughout the unit circle [**GR 995 (9.102.2)**], $E(u)$ is well defined in the integral. While the integral in Eq. (2.16) does not appear in standard references in the literature, it can nevertheless be evaluated since the ellipsoid surface area is known and given by Eq. (1.11). Using Eqs. (1.11) and (1.12), the integral in Eq. (2.16) can be evaluated and expressed as

$$(2.17) \quad \int_0^{\sqrt{\frac{e_1^2-e_2^2}{1-e_2^2}}} \dfrac{(e_1^2-e_2^2)^{\frac{3}{2}} u E(u) du}{\left(e_1^2-e_2^2+e_2^2 u^2\right)^2 \sqrt{(e_1^2-e_2^2)-(1-e_2^2)u^2}}$$

$$= \dfrac{\pi}{4}\left[\sqrt{(1-e_1^2)(1-e_2^2)} + \dfrac{e_1^2 E(\lambda,k)+(1-e_1^2)F(\lambda,k)}{e_1}\right],$$

where $c^2 = ab\sqrt{(1-e_1^2)(1-e_2^2)}$, and

$$(2.18) \quad \lambda = \arcsin e_1, \text{ and } k = \dfrac{e_2}{e_1}, \text{ and } 1 > e_1 > e_2 > 0.$$

The integral above can be characterized by two independent parameters. In Eq. (2.17), those parameters are $(e_1, e_2)$. Instead choose the parameters $(\alpha, k)$ where

(2.19) $\alpha = \sqrt{\dfrac{e_1^2 - e_2^2}{1 - e_2^2}}$ and $k = \dfrac{e_2}{e_1}$,

which implies that,

(2.20) $e_1 = \dfrac{\alpha}{\sqrt{k'^2 + k^2 \alpha^2}}$ and $e_2 = \dfrac{k\alpha}{\sqrt{k'^2 + k^2 \alpha^2}}$, where $k' = \sqrt{1 - k^2}$.

With these definitions, the integral above can be rewritten, so that after a little algebra, Eq. (2.17) becomes

(2.21) $\displaystyle\int_0^\alpha \dfrac{u E(u)\, du}{\left(k'^2 + k^2 u^2\right)^2 \sqrt{\alpha^2 - u^2}} = \dfrac{\pi}{4}\left[\dfrac{\alpha\sqrt{1-\alpha^2}}{\left(k'^2 + k^2 \alpha^2\right)^2} + \dfrac{\alpha^2 E(\lambda, k)}{k'^2\left(k'^2 + k^2 \alpha^2\right)^{\tfrac{3}{2}}} + \dfrac{\left(1 - \alpha^2\right) F(\lambda, k)}{\left(k'^2 + k^2 \alpha^2\right)^{\tfrac{3}{2}}}\right]$,

where $(\lambda, k, k', \alpha)$ are given by,

(2.22) $\lambda = \arcsin\left(\dfrac{\alpha}{\sqrt{k'^2 + k^2 \alpha^2}}\right)$ and $k' = \sqrt{1 - k^2}$, and $1 > k > 0$ and $1 > \alpha > 0$.

The closest published integral relative to that in Eq. (2.21) appears to be due to Prudnikov, Brychkov, and Marichev[27-29] and is rewritten below with parameter changes to better facilitate a comparison with Eq. (2.21), namely,

(2.23) $\displaystyle\int_0^\alpha \dfrac{u E\!\left(\tfrac{u}{\alpha}\right) du}{\left(z^2 + u^2\right)\sqrt{\alpha^2 - u^2}} = \dfrac{\pi \alpha}{2\left(z^2 + \alpha^2\right)} D\!\left(\dfrac{\alpha}{\sqrt{z^2 + \alpha^2}}\right)$, **PR3 185 (2.16.5.5)**

where

(2.24) $z^2 = \dfrac{k'^2}{k^2}$ and $(\alpha, \operatorname{Re} z) > 0$,

and $\boldsymbol{D}(k) \equiv \dfrac{1}{k^2}[\boldsymbol{K}(k) - \boldsymbol{E}(k)]$, where $\boldsymbol{K}(k) \equiv F\!\left(\dfrac{\pi}{2}, k\right)$ is a complete elliptic integrals of the first kind. In general,

(2.25) $D(\lambda, k) = \dfrac{1}{k^2}[F(\lambda, k) - E(\lambda, k)]$, **GR 852 (8.111.5)**

so that, Eq. (2.21) can be equally written as a function of $D(\lambda, k)$ and $F(\lambda, k)$, namely,

$$\text{(2.26)} \quad \int_0^\alpha \frac{u E(u)\,du}{\left(k'^2 + k^2 u^2\right)^2 \sqrt{\alpha^2 - u^2}} = \frac{\pi}{4}\left[\frac{\alpha\sqrt{1-\alpha^2}}{\left(k'^2 + k^2\alpha^2\right)^2} - \frac{\alpha^2 k^2 D(\lambda, k)}{k'^2 \left(k'^2 + k^2\alpha^2\right)^{\frac{3}{2}}} + \frac{F(\lambda, k)}{k'^2 \sqrt{k'^2 + k^2\alpha^2}}\right].$$

The numerical equivalence of integrals found in Eqs. (2.17), (2.21), and (2.26) can be checked by numerical integration and evaluation of the functions. In this study, the numerical equivalence of integrals or series representations was tested by an extensive application of Mathematica $5^{26}$ to the equation in question.

## 2.2  The Integral $I_2(a,b,c) = \int_\alpha^\beta \ln\left(\frac{\varepsilon+u}{\varepsilon-u}\right) \frac{u^2\,du}{\sqrt{(u^2-\alpha^2)(\beta^2-u^2)}}$

By integrating the double integral in Eq. (2.9) with respect to $\theta$ first, the expression for the surface area can be reduced to,

$$\text{(2.27)} \quad S(a,b,c) = 8ab \int_0^{\pi/2} \left[-\frac{\Delta\cos\theta}{2} - \frac{q'^2}{2q} \ln(q\cos\theta + \Delta)\right]_0^{\pi/2} d\phi, \qquad \textbf{GR 182 (2.583.2)}$$

or, equivalently,

$$\text{(2.28)} \quad S(a,b,c) = 2\pi ab + 2ab \int_0^{\pi/2} \ln\left(\frac{1+q}{1-q}\right) \frac{(1-q^2)\,d\phi}{q},$$

where

$$\text{(2.29)} \quad \Delta = \sqrt{1 - q^2 \sin^2\theta} \text{ and } q' = \sqrt{1-q^2},$$

and $q = q(\phi)$ is given by the positive square root of Eq. (2.10). If $q = q(\phi)$ is inverted and $\phi$ differentiated with respect to $q$, then Eq. (2.28) can be rewritten in terms of $q$ to become

$$\text{(2.30)} \quad S(a,b,c) = 2\pi ab + 2ab \int_{e_2}^{e_1} \ln\left(\frac{1+q}{1-q}\right) \frac{(1-q^2)\,dq}{\sqrt{(e_1^2 - q^2)(q^2 - e_2^2)}},$$

where

$$\text{(2.31)} \quad \phi = \arcsin\sqrt{\frac{e_1^2 - q^2}{e_1^2 - e_2^2}} \text{ and } 1 > e_1 > e_2 > 0$$

have been applied. Only one of the two integral parts in Eq. (2.30) could be found in the literature; namely,

$$\text{(2.32)} \quad \int_{e_2}^{e_1} \ln\left(\frac{1+q}{1-q}\right) \frac{dq}{\sqrt{(e_1^2-q^2)(q^2-e_2^2)}} = \frac{\pi}{e_1} F(\varphi,k), \qquad \textbf{PR1 510 (2.6.13.2)}$$

where

$$\text{(2.33)} \quad \varphi = \arcsin e_1 \text{ and } k = \frac{e_2}{e_1}.$$

With some algebra, the remaining integral part involving $q^2$ can be evaluated from Eqs. (1.11) and (2.32), and is given below,

$$\text{(2.34)} \quad \int_{e_2}^{e_1} \ln\left(\frac{1+q}{1-q}\right) \frac{q^2 dq}{\sqrt{(e_1^2-q^2)(q^2-e_2^2)}} = \pi\left[1-\sqrt{(1-e_1^2)(1-e_2^2)}\right] + \pi e_1\left[F(\varphi,k) - E(\varphi,k)\right].$$

### 2.2.1 The Series $\Sigma_1(a,b,c) = \int_{e_2}^{e_1} \ln\left(\frac{1+q}{1-q}\right) \frac{dq}{\sqrt{(e_1^2-q^2)(q^2-e_2^2)}}$

The integrals found in Eqs. (2.32) and (2.34) are interesting in that an infinite series representation can be readily found by expanding the natural logarithm about $q = 0$ and then applying a recurrence relation to the subsequent series of integrals. Let

$$\text{(2.35)} \quad \ln\left(\frac{1+q}{1-q}\right) = 2\sum_{m=0}^{\infty} \frac{q^{2m+1}}{2m+1}, \qquad \textbf{GR 51 (1.513.1)}$$

where $q^2 < 1$ from Eq. (2.10). Inserting Eq. (2.35) into Eq. (2.32) and interchanging the sum and integral allows Eq. (2.32) to be written as an infinite series of integrals,

$$\text{(2.36)} \quad \int_{e_2}^{e_1} \ln\left(\frac{1+q}{1-q}\right) \frac{dq}{\sqrt{(e_1^2-q^2)(q^2-e_2^2)}} = \sum_{m=0}^{\infty} \frac{2}{(2m+1)} \int_{e_2}^{e_1} \frac{q^{2m+1} dq}{\sqrt{(e_1^2-q^2)(q^2-e_2^2)}},$$

where

$$\text{(2.37)} \quad \int_{e_2}^{e_1} \frac{q^{2m+1} dq}{\sqrt{(e_1^2-q^2)(q^2-e_2^2)}} = \frac{e_2^{2m+1}}{e_1} \int_{1}^{1/k} \frac{t^{2m+1} dt}{\sqrt{(t^2-1)(1-k^2 t^2)}},$$

and

$$\text{(2.38)} \quad t = \frac{q}{e_2} \text{ and } k = \frac{e_2}{e_1}.$$

The last integral in Eq. (2.37) can be evaluated and leads to the recurrence relation below,

(2.39) $\quad I_{2m+1} = \int \frac{t^{2m+1} dt}{\sqrt{(t^2-1)(1-k^2 t^2)}} = \int \frac{d\phi}{(1-k'^2 \sin^2 \phi)^{m+1}} = \int \mathrm{nd}^{2m+1} u\, du,$  **BY 194**

where $t = \dfrac{1}{\sqrt{1-k'^2 \sin^2 \phi}} = \dfrac{1}{\sqrt{1-k'^2 \mathrm{sn}^2(u,k')}} = \dfrac{1}{\mathrm{dn}(u,k')} \equiv \mathrm{nd}\, u$, and  **BY 194**

(2.40) $\quad I_1 = \dfrac{1}{k} \arctan\left( \dfrac{k\, \mathrm{sn}\, u - \mathrm{cn}\, u}{k\, \mathrm{sn}\, u + \mathrm{cn}\, u} \right),$  **BY 194 (315.01)**

(2.41) $\quad I_3 = \dfrac{1}{2k^3}\left\{ (2-k'^2)\arctan\left( \dfrac{k\, \mathrm{sn}\, u - \mathrm{cn}\, u}{k\, \mathrm{sn}\, u + \mathrm{cn}\, u} \right) - kk'^2 \, \mathrm{sn}\, u\, \mathrm{cn}\, u\, \mathrm{nd}^2 u \right\},$  **BY 194 (315.03)**

(2.42) $\quad I_{2n+5} = \dfrac{(2n+3)(2-k'^2) I_{2n+3} - 2(n+1) I_{2n+1} - k'^2\, \mathrm{sn}\, u\, \mathrm{cn}\, u\, \mathrm{nd}^{2n+4} u}{2(n+2)k^2},$  **BY 195 (315.06)**

where $n = 0,1,2,3,\ldots$, $m = n+2$, and $\mathrm{sn}\, u$, $\mathrm{cn}\, u$, and $\mathrm{nd}\, u$ are elliptic functions. Note that relative to **BY 194 – 195**, we have exchanged $k \leftrightarrow k'$ to be consistent with notation used here, namely $k = \dfrac{e_2}{e_1}$, set $m \to 2m+1$ in Eq. (2.39), so that only positive, odd powers of $t$ are considered, and set $n \to n+1$ in Eq. (2.42), so that all $I$ terms are defined and consistent when $m = 0,1,2,3,\ldots$.

Several simplifications are possible at this point. These simplifications are best revealed by expressing the elliptic functions in terms of $q$ rather than $t$, and then considering the upper and lower limits of the $q$ integral in Eq. (2.37), namely $(e_1, e_2)$. With these modifications and the usual algebra of elliptic functions, we find that

(2.43) $\quad \mathrm{nd}\, u = \dfrac{q}{e_2}, \quad \mathrm{sn}\, u = \dfrac{e_1}{q}\sqrt{\dfrac{q^2 - e_2^2}{e_1^2 - e_2^2}}, \quad \text{and} \quad \mathrm{cn}\, u = \dfrac{e_2}{q}\sqrt{\dfrac{e_1^2 - q^2}{e_1^2 - e_2^2}},$

(2.44) $\quad \mathrm{sn}\, u\, \mathrm{cn}\, u\, \mathrm{nd}^{2n} u = \dfrac{e_1 e_2}{q^2}\sqrt{\dfrac{e_1^2 - q^2}{e_1^2 - e_2^2}}\sqrt{\dfrac{q^2 - e_2^2}{e_1^2 - e_2^2}}\left(\dfrac{q}{e_2}\right)^{2n},$

and,

(2.45) $\quad \dfrac{k\, \mathrm{sn}\, u - \mathrm{cn}\, u}{k\, \mathrm{sn}\, u + \mathrm{cn}\, u} = \dfrac{\sqrt{q^2 - e_2^2} - \sqrt{e_1^2 - q^2}}{\sqrt{q^2 - e_2^2} + \sqrt{e_1^2 - q^2}}.$

From Eq. (2.44), it is clear that $\operatorname{sn} u \operatorname{cn} u \operatorname{nd}^{2n} u$ vanishes on the upper and lower limits $(e_1, e_2)$, so that these terms may be dropped from further consideration. From Eq. (2.45), it is clear that

$$(2.46) \quad \arctan\left(\frac{k \operatorname{sn} u - \operatorname{cn} u}{k \operatorname{sn} u + \operatorname{cn} u}\right) = \begin{Bmatrix} +\dfrac{\pi}{4} \\ -\dfrac{\pi}{4} \end{Bmatrix}, \text{ when } q = \begin{Bmatrix} e_1 \\ e_2 \end{Bmatrix},$$

so that $\arctan\left(\dfrac{k \operatorname{sn} u - \operatorname{cn} u}{k \operatorname{sn} u + \operatorname{cn} u}\right)$ contributes $\dfrac{\pi}{2}$ every time it is evaluated as the difference between the upper and lower limits $(e_1, e_2)$. Furthermore, note that for each $\tilde{I}_{2m+1} \equiv I_{2m+1}(e_1) - I_{2m+1}(e_2)$,

$$(2.47) \quad \tilde{I}_1 = \frac{e_1}{e_2}\left(\frac{\pi}{2}\right),$$

$$(2.48) \quad \tilde{I}_3 = \frac{e_1(e_1^2 + e_2^2)}{2 e_2^3}\left(\frac{\pi}{2}\right),$$

$$(2.49) \quad \tilde{I}_{2n+5} = \frac{(e_1^2 + e_2^2)(2n+3)\tilde{I}_{2n+3} - 2 e_1^2(n+1)\tilde{I}_{2n+1}}{2(n+2) e_2^2}, \quad n = 0,1,2,3,\ldots, \quad m = n+2,$$

contains a factor of $\dfrac{e_1}{e_2^{2m+1}}$ that exactly cancels the factor of $\dfrac{e_2^{2m+1}}{e_1}$ that appears in Eq. (2.37), and that each $\tilde{I}_{2m+1}$ contains a factor of $\dfrac{\pi}{2}$. Consequently, after multiplying each term in the series by $\dfrac{2}{(2m+1)}$, the integral in Eq. (2.36) can be represented by

$$(2.50) \quad \int_{e_2}^{e_1} \ln\left(\frac{1+q}{1-q}\right) \frac{dq}{\sqrt{(e_1^2 - q^2)(q^2 - e_2^2)}} = \pi \sum_{m=0}^{\infty} \Lambda_{2m+1},$$

where

$$(2.51) \quad \Lambda_1 = 1,$$

$$(2.52) \quad \Lambda_3 = \frac{e_1^2 + e_2^2}{6},$$

$$(2.53) \quad \Lambda_{2n+5} = \frac{(e_1^2 + e_2^2)(2n+3)^2 \Lambda_{2n+3} - 2 e_1^2 e_2^2 (n+1)(2n+1)\Lambda_{2n+1}}{2(n+2)(2n+5)}, \quad n = 0,1,2,3,\ldots, \quad m = n+2,$$

and $\Lambda_{2n+1} = \dfrac{\tilde{I}_{2n+1}}{(2n+1)}$ but $e_2^2$ multiplies only $\Lambda_{2n+1}$ in Eq. (2.53) as part of the algebra that removes $e_2^2$ from the denominator in Eq. (2.49).

Since $F(\arcsin e_1, k)$ is an odd function of $\arcsin e_1$ [**GR 854 (8.121.1)**], it can be expanded as a Maclaurin (Taylor) series of odd powers of $e_1$ about the point $e_1 = 0$ in Eq. (2.32). Let

(2.54) $\quad F(\arcsin e_1, k) = \displaystyle\sum_{m=0}^{\infty} \dfrac{e_1^{2m+1} F^{2m+1}(\arcsin 0, k)}{(2m+1)!}$,

where $F^{2m+1}(\arcsin 0, k)$ represents $\displaystyle\lim_{z \to 0}\left[\dfrac{\partial^{2m+1} F(\arcsin z, k)}{\partial z^{2m+1}}\right]$, then Eqs. (2.32), (2.33), and (2.54) imply a term-by-term equivalence with the infinite series defined by Eqs. (2.50) – (2.53), so that the $F^{2m+1}(\arcsin 0, k)$ can be expressed as

(2.55) $\quad F^{2m+1}(\arcsin 0, k) = \dfrac{(2m+1)!\,\Lambda_{2m+1}}{e_1^{2m}}$, for $m = 0, 1, 2, 3, \ldots$,

where

$k = \dfrac{e_2}{e_1}$, and the $\Lambda_{2m+1}$ are given by Eqs. (2.51) – (2.53).

### 2.2.2 The Series $\Sigma_2(a,b,c) = \displaystyle\int_{e_2}^{e_1} \ln\left(\dfrac{1+q}{1-q}\right)\dfrac{q^2\, dq}{\sqrt{(e_1^2 - q^2)(q^2 - e_2^2)}}$

The derivation of the series $\Sigma_2(a,b,c)$ quickly follows from $\Sigma_1(a,b,c)$. Here Eqs. (2.36) – (2.38) become

(2.56) $\quad \displaystyle\int_{e_2}^{e_1} \ln\left(\dfrac{1+q}{1-q}\right)\dfrac{q^2\, dq}{\sqrt{(e_1^2 - q^2)(q^2 - e_2^2)}} = \sum_{m=0}^{\infty} \dfrac{2}{(2m+1)}\int_{e_2}^{e_1} \dfrac{q^{2m+3}\, dq}{\sqrt{(e_1^2 - q^2)(q^2 - e_2^2)}}$,

where

(2.57) $\quad \displaystyle\int_{e_2}^{e_1} \dfrac{q^{2m+3}\, dq}{\sqrt{(e_1^2 - q^2)(q^2 - e_2^2)}} = \dfrac{e_2^{2m+3}}{e_1}\int_{1}^{1/k} \dfrac{t^{2m+3}\, dt}{\sqrt{(t^2-1)(1-k^2 t^2)}}$,

and

(2.58) $t = \dfrac{q}{e_2}$ and $k = \dfrac{e_2}{e_1}$.

Again, the last integral in Eq. (2.57) can be evaluated, and leads – with the same simplifications – to the same recurrence relation given by Eqs. (2.47) – (2.49). At this point, the $q^{2m+3}$ in Eq. (2.56) are multiplied by $\dfrac{2}{(2m+1)}$ rather than $\dfrac{2}{(2m+3)}$. This small modification produces a slightly different series, but a substantially different functional representation for the integral. Let $\Omega_{2n+1} = \dfrac{\tilde{I}_{2n+1}}{|2n-1|}$ in Eqs. (2.47) – (2.49), then with the remaining simplifications, the integral in Eq. (2.56) can be represented by

(2.59) $\displaystyle\int_{e_2}^{e_1} \ln\left(\dfrac{1+q}{1-q}\right) \dfrac{q^2 \, dq}{\sqrt{(e_1^2 - q^2)(q^2 - e_2^2)}} = \pi \left(\sum_{m=0}^{\infty} \Omega_{2m+1} - 1\right)$,

where

(2.60) $\Omega_1 = 1$,

(2.61) $\Omega_3 = \dfrac{e_1^2 + e_2^2}{2}$,

(2.62) $\Omega_{2n+5} = \dfrac{(e_1^2 + e_2^2)(2n+1)(2n+3)\Omega_{2n+3} - 2e_1^2 e_2^2 (n+1)|2n-1|\Omega_{2n+1}}{2(n+2)(2n+3)}$, $n = 0,1,2,3,\ldots$, $m = n+2$,

and, as with Equation (2.53), $e_2^2$ multiplies only $\Omega_{2n+1}$ in Eq. (2.62) as part of the algebra that removes $e_2^2$ from the denominator in Eq. (2.49).

It is instructive to compare the first few terms of the series given by Eqs. (2.59) – (2.62) with the corresponding series derived by expanding the functions in Eq. (2.34) in powers of $e_1$ and $e_2$. By dropping the common factor $\pi$, and rewriting Eq. (2.34) in terms of the series from Eqs. (2.59) – (2.62), one has the expression,

(2.63) $\left(\displaystyle\sum_{m=0}^{\infty} \Omega_{2m+1} - 1\right) = \left(1 - \sqrt{1-e_1^2}\sqrt{1-e_2^2}\right) + e_1\left[F(\varphi,k) - E(\varphi,k)\right]$,

where the first few non-zero terms from $\left(\displaystyle\sum_{m=0}^{\infty} \Omega_{2m+1} - 1\right)$ are

(2.64) $\Omega_3 = \dfrac{e_1^2 + e_2^2}{2}$,

$$(2.65) \quad \Omega_5 = \frac{3e_1^4 + 2e_1^2 e_2^2 + 3e_2^4}{24},$$

$$(2.66) \quad \Omega_7 = \frac{5e_1^6 + 3e_1^4 e_2^2 + 3e_1^2 e_2^4 + 5e_2^6}{80}.$$

With the $\sqrt{1-e_1^2}$ and $\sqrt{1-e_2^2}$ square root terms expressed as Maclaurin series in $e_1$ and $e_2$, respectively, and the product of the two infinite series arranged as a Cauchy product, the first few non-zero terms of $\left(1 - \sqrt{1-e_1^2}\sqrt{1-e_2^2}\right) = \sum_{m=1}^{\infty} \Theta_{2m+1}$ can be written as

$$(2.67) \quad \Theta_3 = \frac{e_1^2 + e_2^2}{2},$$

$$(2.68) \quad \Theta_5 = \frac{e_1^4 - 2e_1^2 e_2^2 + e_2^4}{8},$$

$$(2.69) \quad \Theta_7 = \frac{e_1^6 - e_1^4 e_2^2 - e_1^2 e_2^4 + e_2^6}{16},$$

while the first few non-zero terms of $e_1\left[F(\arcsin e_1, k) - E(\arcsin e_1, k)\right] = \sum_{m=2}^{\infty} \Psi_{2m+1}$, expressed as a Maclaurin series in $e_1$, can be written as

$$(2.70) \quad \Psi_5 = \frac{e_1^2 e_2^2}{3},$$

$$(2.71) \quad \Psi_7 = \frac{e_1^4 e_2^2 + e_1^2 e_2^4}{10},$$

where $k = \frac{e_2}{e_1}$. Note that there is a clear correspondence between $e_1^{n_1} e_2^{n_2}$ terms of equal power in $n_1 + n_2$ from Eqs. (2.64) – (2.66) and like terms from the sum of Eqs. (2.67) – (2.69) and (2.70) – (2.71).

Finally, we rewrite Eq. (2.34) in a more general form so that it more closely follows an expression given by Prudnikov, Brychkov, and Marichev[27-29], namely,

$$(2.72) \quad \int_\alpha^\beta \ln\left(\frac{\varepsilon+u}{\varepsilon-u}\right) \frac{du}{\sqrt{(u^2 - \alpha^2)(\beta^2 - u^2)}} = \frac{\pi}{\beta} F\left(\arcsin\frac{\beta}{\varepsilon}, \frac{\alpha}{\beta}\right), \quad \textbf{PR1 510 (2.6.13.2)}$$

where $|\varepsilon| > |\beta|$ and $\beta > \alpha$. In Eq. (2.32) above, we had set $\alpha = e_2$, $\beta = e_1$, and $\varepsilon = 1$, where $1 > e_1 > e_2 > 0$. If we now let $q = \dfrac{u}{\varepsilon}$, $e_2 = \dfrac{\alpha}{\varepsilon}$, and $e_1 = \dfrac{\beta}{\varepsilon}$, then the integral in Eq. (2.34) can be rewritten as

$$(2.73) \quad \int_{e_2}^{e_1} \ln\left(\frac{1+q}{1-q}\right) \frac{q^2\,dq}{\sqrt{(e_1^2 - q^2)(q^2 - e_2^2)}} = \frac{1}{\varepsilon}\int_{\alpha}^{\beta} \ln\left(\frac{\varepsilon+u}{\varepsilon-u}\right) \frac{u^2\,du}{\sqrt{(u^2 - \alpha^2)(\beta^2 - u^2)}}.$$

Multiplying through by $\varepsilon$ and rewriting the remaining functions leads to a more general form for Eq. (2.34), namely,

$$(2.74) \quad \int_{\alpha}^{\beta} \ln\left(\frac{\varepsilon+u}{\varepsilon-u}\right) \frac{u^2\,du}{\sqrt{(u^2 - \alpha^2)(\beta^2 - u^2)}} = \frac{\pi}{\varepsilon}\left(\varepsilon^2 - \sqrt{(\varepsilon^2 - \alpha^2)(\varepsilon^2 - \beta^2)}\right) + \pi\beta\left[F(\varphi,k) - E(\varphi,k)\right],$$

$$(2.75) \quad \varphi = \arcsin\frac{\beta}{\varepsilon} \text{ and } k = \frac{\alpha}{\beta},$$

where $1 > e_1 > e_2 > 0$ implies that $\varepsilon > \beta > \alpha > 0$.

### 2.3 The Integral $I_3(a,b,c) = \displaystyle\int_0^{\pi/2} \dfrac{E(u,k')\sin u \cos u \, du}{\left(1 - k'^2 \cosh^2 v \sin^2 u\right)\sqrt{1 - k'^2 \sin^2 u}}$

In the last section, the integral in Eq. (2.34) was evaluated using the known surface area of an ellipsoid where $a > b > c > 0$. To find related integrals, rewrite the integral part of Eq. (2.34) such that,

$$(2.76) \quad \int_{e_2}^{e_1} \ln\left(\frac{1+q}{1-q}\right) \frac{q^2\,dq}{\sqrt{(e_1^2 - q^2)(q^2 - e_2^2)}} = 2\int_{e_2}^{e_1} \frac{q^2\, \text{arctanh}\, q \, dq}{\sqrt{(e_1^2 - q^2)(q^2 - e_2^2)}}, \quad \text{GR 54 (1.622.7)}$$

and then integrate by parts, where

$$(2.77) \quad \omega = \text{arctanh}\, q,$$

and

$$(2.78) \quad d\xi = \frac{q^2\,dq}{\sqrt{(e_1^2 - q^2)(q^2 - e_2^2)}}.$$

Since the integration of Eq. (2.78) involves an elliptic integral of the second kind, two functionally different but equivalent solutions for $\xi$ can be defined on the interval $[e_2, e_1]$ depending upon whether the upper limit or lower limit of the integral is taken as a variable. Let $\xi_i$, where $i = 1, 2$, represent those two solutions for $\xi$, so that with the upper limit $q$ as the variable, the first solution for $\xi$ can be defined by

$$(2.79) \quad \xi_1 = \begin{cases} \displaystyle\int_{e_2}^{q} \frac{\tilde{q}^2 \, d\tilde{q}}{\sqrt{(e_1^2 - \tilde{q}^2)(\tilde{q}^2 - e_2^2)}}, \\ 0 \end{cases}$$

$$= \begin{cases} e_1 E(u, k') - \dfrac{1}{q}\sqrt{(e_1^2 - q^2)(q^2 - e_2^2)}, & \text{when} \quad e_1 \geq q > e_2 > 0, \\ 0, & \text{when} \quad q = e_2, \end{cases} \qquad \textbf{GR 277 (3.153.7)}$$

where

$$(2.80) \quad u = \arcsin\left(\frac{e_1}{q}\sqrt{\frac{q^2 - e_2^2}{e_1^2 - e_2^2}}\right) \quad \text{and} \quad k' = \sqrt{1 - \frac{e_2^2}{e_1^2}}.$$

By applying the standard technique for integration by parts, namely,

$$(2.81) \quad \int_{e_2}^{e_1} \omega \, d\xi_i = [\omega \xi_i]_{e_2}^{e_1} - \int_{e_2}^{e_1} \xi_i \, d\omega,$$

Eq. (2.76) reduces to

(2.82)

$$\int_{e_2}^{e_1} \ln\left(\frac{1+q}{1-q}\right) \frac{q^2 \, dq}{\sqrt{(e_1^2 - q^2)(q^2 - e_2^2)}} = 2e_1 \mathbf{E}(k') \operatorname{arctanh} e_1 - 2e_1 \int_{e_2}^{e_1} \frac{E(u, k') \, dq}{1 - q^2} + 2 \int_{e_2}^{e_1} \frac{\sqrt{(e_1^2 - q^2)(q^2 - e_2^2)} \, dq}{q(1 - q^2)},$$

where Eqs. (2.77) and (2.79) define $\omega$ and $\xi_1$, respectively. A second solution for $\xi$ denoted by $\xi_2$ and satisfying Eq. (2.78) is discussed in Section 2.4. However no other solutions for $\xi$, for example, $\xi_3 = \lambda \xi_1 + (1 - \lambda) \xi_2$, where $1 > \lambda > 0$, are considered in this study.

The last integral in Eq. (2.82) is a pseudo-elliptic integral that can be readily integrated to elementary functions. Consider the last integral first. Multiplying numerator and denominator by $\sqrt{(e_1^2 - q^2)(q^2 - e_2^2)}$, removing the whole part, and writing the remainder in terms of partial fractions, reduces the pseudo-elliptic integral to a form that can be readily integrated, namely,

$$(2.83) \quad \int_{e_2}^{e_1} \frac{\sqrt{(e_1^2-q^2)(q^2-e_2^2)}\,dq}{q(1-q^2)} = \int_{e_2}^{e_1} \frac{(e_1^2-q^2)(q^2-e_2^2)\,dq}{q(1-q^2)\sqrt{(e_1^2-q^2)(q^2-e_2^2)}},$$

where

$$(2.84) \quad \frac{(e_1^2-q^2)(q^2-e_2^2)}{q(1-q^2)} = q - \frac{e_1^2 e_2^2}{q} - \frac{(1-e_1^2)(1-e_2^2)}{2(1-q)} + \frac{(1-e_1^2)(1-e_2^2)}{2(1+q)},$$

and

$$(2.85) \quad \int_{e_2}^{e_1} \frac{q\,dq}{\sqrt{(e_1^2-q^2)(q^2-e_2^2)}} = \frac{1}{2}\int_{e_2^2}^{e_1^2} \frac{d\tau}{\sqrt{(e_1^2-\tau)(\tau-e_2^2)}} = \frac{\pi}{2}, \qquad \textbf{GR 93 (2.261)}$$

$$(2.86) \quad \int_{e_2}^{e_1} \frac{-e_1^2 e_2^2\,dq}{q\sqrt{(e_1^2-q^2)(q^2-e_2^2)}} = \frac{-e_1^2 e_2^2}{2}\int_{e_2^2}^{e_1^2} \frac{d\tau}{\tau\sqrt{(e_1^2-\tau)(\tau-e_2^2)}} = \frac{-\pi e_1 e_2}{2}, \qquad \textbf{GR 95 (2.266)}$$

$$(2.87) \quad \int_{e_2}^{e_1} \frac{-(1-e_1^2)(1-e_2^2)\,q\,dq}{(1-q^2)\sqrt{(e_1^2-q^2)(q^2-e_2^2)}} = \frac{-(1-e_1^2)(1-e_2^2)}{2}\int_{1-e_1^2}^{1-e_2^2} \frac{d\eta}{\eta\sqrt{[(1-e_2^2)-\eta][\eta-(1-e_1^2)]}}$$

$$= \frac{-\pi\sqrt{(1-e_1^2)(1-e_2^2)}}{2}, \qquad \textbf{GR 95 (2.266)}$$

where $\tau = q^2$ has been used in Eqs. (2.85) and (2.86), and the last two terms of Eq. (2.84) have been combined and $\eta = 1-q^2$ has been applied to Eq. (2.87). Consequently we have from Eqs. (2.83) and (2.84),

$$(2.88) \quad \int_{e_2}^{e_1} \frac{\sqrt{(e_1^2-q^2)(q^2-e_2^2)}\,dq}{q(1-q^2)} = \frac{\pi}{2}\left[1 - e_1 e_2 - \sqrt{(1-e_1^2)(1-e_2^2)}\right],$$

so that with Eq. (2.34), the integral containing $E(u,k')$ in Eq. (2.82) can be isolated and expressed as

$$(2.89) \quad \int_{e_2}^{e_1} \frac{E(u,k')\,dq}{1-q^2} = \boldsymbol{E}(k')\operatorname{arctanh} e_1 - \frac{\pi e_2}{2} - \frac{\pi}{2}\left[F(\varphi,k) - E(\varphi,k)\right],$$

where

$$(2.90) \quad u = \arcsin\left(\frac{e_1}{q}\sqrt{\frac{q^2 - e_2^2}{e_1^2 - e_2^2}}\right), \quad k' = \sqrt{1-k^2}, \quad \varphi = \arcsin e_1, \text{ and } k = \frac{e_2}{e_1}.$$

The integral in Eq. (2.89) can be written in a more conventional form by expressing $q$ as a function of $u$, namely,

$$(2.91) \quad q = \frac{e_1 e_2}{\sqrt{e_1^2 - (e_1^2 - e_2^2)\sin^2 u}},$$

and rewriting Eq. (2.89) in terms of $u$, so that,

$$(2.92) \quad \int_{e_2}^{e_1} \frac{E(u,k')\, dq}{1 - q^2} = \frac{e_2 k'^2}{(1 - e_2^2)} \int_0^{\pi/2} \frac{E(u,k') \sin u \cos u\, du}{\left[1 - k'^2\left(\frac{1}{1-e_2^2}\right)\sin^2 u\right]\sqrt{1 - k'^2 \sin^2 u}}.$$

The last integral in Eq. (2.92) can be characterized in a number of ways. Let

$$(2.93) \quad \cosh v = \frac{1}{\sqrt{1 - e_2^2}} \quad \text{or, equivalently,} \quad \tanh v = e_2,$$

then the integral in Eq. (2.92) can be characterized by either $(v,k)$ or $(v,k')$. With Eqs. (2.92) and (2.93) applied to Eq. (2.89), the integral can be expressed as

$$(2.94) \quad \int_0^{\pi/2} \frac{E(u,k')\sin u \cos u\, du}{(1 - k'^2 \cosh^2 v \sin^2 u)\sqrt{1 - k'^2 \sin^2 u}}$$

$$= \frac{1}{k'^2 \sinh v \cosh v}\left\{E(k')\operatorname{arctanh}\left(\frac{\tanh v}{k}\right) - \frac{\pi \tanh v}{2} - \frac{\pi}{2}[F(\varphi,k) - E(\varphi,k)]\right\},$$

where

$$(2.95) \quad \varphi = \arcsin\left(\frac{\tanh v}{k}\right) \quad \text{and} \quad k' = \sqrt{1 - k^2},$$

and $1 > e_1 > e_2 > 0$ implies that $1 > k > \tanh v > 0$.

## 2.4 The Integral $I_4(a,b,c) = \int_0^{\pi/2} \dfrac{E(u,k')\sin u \cos u \, du}{(1+k'^2 \sinh^2\mu \sin^2 u)\sqrt{1-k'^2 \sin^2 u}}$

An integral related to that in Eq. (2.94) could be obtained by integrating Eq. (2.76) by parts but with the lower limit as the variable in solving Eq. (2.78) for $\xi$. Equations (2.76) – (2.78) are repeated for reference, so that

(2.96) $\displaystyle\int_{e_2}^{e_1} \ln\left(\dfrac{1+q}{1-q}\right) \dfrac{q^2 \, dq}{\sqrt{(e_1^2 - q^2)(q^2 - e_2^2)}} = 2\int_{e_2}^{e_1} \dfrac{q^2 \operatorname{arctanh} q \, dq}{\sqrt{(e_1^2 - q^2)(q^2 - e_2^2)}}$,      **GR 54 (1.622.7)**

where

(2.97) $\omega = \operatorname{arctanh} q$,

(2.98) $d\xi = \dfrac{q^2 \, dq}{\sqrt{(e_1^2 - q^2)(q^2 - e_2^2)}}$.

With the lower limit $q$ as the variable, the second solution for $\xi$ can be defined by,

(2.99)

$-\xi_2 = \begin{cases} 0, & \text{when} \quad q = e_1 \\ \displaystyle\int_q^{e_1} \dfrac{\tilde{q}^2 \, d\tilde{q}}{\sqrt{(e_1^2 - \tilde{q}^2)(\tilde{q}^2 - e_2^2)}} = e_1 E(u,k'), & \text{when} \quad e_1 > q \geq e_2 > 0, \end{cases}$      **GR 277 (3.153.8)**

where

(2.100) $u = \arcsin \sqrt{\dfrac{e_1^2 - q^2}{e_1^2 - e_2^2}}$ and $k' = \sqrt{1 - \dfrac{e_2^2}{e_1^2}}$.

The minus sign in Eq. (2.99) arises because the integration interval is reversed, so that the partial derivative $-e_1 \dfrac{\partial E(u,k')}{\partial q}$ equals the integrand when $e_1 > q \geq e_2 > 0$. The two solutions for $\xi$, $\xi_1$ and $\xi_2$, given by Eqs. (2.79) and (2.99), respectively, are related to one another by the addition formula for $E(\varphi_i, k')$, where $i = 1, 2$, namely,

(2.101) $E(\varphi_1, k') \pm E(k') = E(\varphi_2, k') \pm k'^2 \sin\varphi_1 \sin\varphi_2$,      **BY 13 (166.01)**

where

$$(2.102) \cos\varphi_2 = \mp \frac{\sin\varphi_1 \sqrt{(1-k'^2 \sin^2\varphi_1)(1-k'^2)}}{1-k'^2 \sin^2\varphi_1}. \qquad \text{BY 13 (166.01)}$$

Using Eqs. (2.80) and (2.100) for $\varphi_1$ and $\varphi_2$, respectively, Eq. (2.101) becomes

$$(2.103)\ e_1 E\left(\arcsin\left(\frac{e_1}{q}\sqrt{\frac{q^2 - e_2^2}{e_1^2 - e_2^2}}\right), k'\right) - \frac{1}{q}\sqrt{(e_1^2 - q^2)(q^2 - e_2^2)} = e_1 E\left(\arcsin\sqrt{\frac{e_1^2 - q^2}{e_1^2 - e_2^2}}, k'\right) - e_1 \boldsymbol{E}(k'),$$

so that the two solutions for $\xi$ that satisfy Eq. (2.78) or (2.98), namely $\xi_1$ and $\xi_2$, given by Eqs. (2.79) and (2.99), are equivalent, within an integration constant, on the interval $[e_2, e_1]$.

Following the procedure outlined in Section 2.3, with integration by parts defined by Eq. (2.81), and $\omega$ and $\xi_2$ defined by Eqs. (2.97) and (2.99), respectively, Eq. (2.96) reduces to

$$(2.104) \int_{e_2}^{e_1} \ln\left(\frac{1+q}{1-q}\right) \frac{q^2 \, dq}{\sqrt{(e_1^2 - q^2)(q^2 - e_2^2)}} = 2e_1 \boldsymbol{E}(k') \operatorname{arctanh} e_2 + 2e_1 \int_{e_2}^{e_1} \frac{E(u, k') \, dq}{1 - q^2},$$

where $(u, k')$ are given by Eq. (2.100). Again, from Eq. (2.34), the integral containing $E(u, k')$ in Eq. (2.104) can be isolated and can be expressed as

$$(2.105) \int_{e_2}^{e_1} \frac{E(u, k') \, dq}{1 - q^2} = \frac{\pi}{2 e_1}\left[1 - \sqrt{(1 - e_1^2)(1 - e_2^2)}\right] + \frac{\pi}{2}\left[F(\varphi, k) - E(\varphi, k)\right] - \boldsymbol{E}(k') \operatorname{arctanh} e_2,$$

where

$$(2.106)\ u = \arcsin\sqrt{\frac{e_1^2 - q^2}{e_1^2 - e_2^2}},\ k' = \sqrt{1 - k^2},\ \varphi = \arcsin e_1,\ \text{and}\ k = \frac{e_2}{e_1}.$$

The integral in Eq. (2.105) can be rewritten in terms of $u$ with

$$(2.107)\ q = e_1 \sqrt{1 - k'^2 \sin^2 u},\ \text{so that,}$$

$$(2.108) \int_{e_2}^{e_1} \frac{E(u, k') \, dq}{1 - q^2} = \frac{e_1 k'^2}{(1 - e_1^2)} \int_0^{\pi/2} \frac{E(u, k') \sin u \cos u \, du}{\left[1 + k'^2\left(\frac{e_1^2}{1 - e_1^2}\right) \sin^2 u\right] \sqrt{1 - k'^2 \sin^2 u}}.$$

The last integral in Eq. (2.108) can be characterized in a number of ways. Let

(2.109) $\sinh \mu = \dfrac{e_1}{\sqrt{1-e_1^2}}$ or, equivalently, $\tanh \mu = e_1$,

then $k \tanh \mu = e_2$, and the integral in Eq. (2.108) can be characterized by either $(\mu, k)$ or $(\mu, k')$. Using Eqs. (2.108) and (2.109) applied to Eq. (2.105), the integral can be expressed as

$$(2.110) \int_0^{\pi/2} \frac{E(u,k') \sin u \cos u \, du}{\left(1 + k'^2 \sinh^2 \mu \sin^2 u\right)\sqrt{1 - k'^2 \sin^2 u}}$$

$$= \frac{-1}{k'^2 \sinh \mu \cosh \mu} \left\{ E(k') \operatorname{arctanh}(k \tanh \mu) - \frac{\pi}{2}\left[F(\varphi,k) - E(\varphi,k) + \tanh \mu \sqrt{1 + k'^2 \sinh^2 \mu}\right] \right.$$

$$\left. - \frac{\pi}{2} \coth \mu \left(1 - \sqrt{1 + k'^2 \sinh^2 \mu}\right) \right\},$$

where

(2.111) $\varphi = \arcsin(\tanh \mu)$ and $k' = \sqrt{1 - k^2}$,

and $1 > e_1 > e_2 > 0$ implies that $1 > \tanh \mu > 0$ and $1 > k > 0$.

Although the right side of Eq. (2.110) can be further simplified, we preserve the present form to facilitate the analytic extension of Eq. (2.110) as described in Section 4.4.

## 2.5  The Integral $I_5(a,b,c) = \displaystyle\int_0^{\pi/2} \frac{F(u,k') \sin u \cos u \, du}{\left(1 + k'^2 \sinh^2 \mu \sin^2 u\right)\sqrt{1 - k'^2 \sin^2 u}}$

As indicated above, there are many ways to manipulate integrals, but one more deserves mention in this study. In July 2003, Dieckmann[19] issued a communication that contained a variable substitution that is useful here, although for a different reason than originally proposed. What is interesting for this study is that Dieckmann's variable substitution, given below, allows for the derivation and evaluation of $F(u,k')$ integrals corresponding to the $E(u,k')$ integrals found in Sections 2.3 and 2.4. With Eq. (2.30) repeated for reference,

$$(2.112) \; S(a,b,c) = 2\pi ab + 2ab \int_{e_2}^{e_1} \ln\left(\frac{1+q}{1-q}\right) \frac{(1-q^2)\,dq}{\sqrt{(e_1^2 - q^2)(q^2 - e_2^2)}},$$

make the variable substitution $w = 1 - q^2$, so that the integral becomes

$$(2.113) \int_{e_2}^{e_1} \ln\left(\frac{1+q}{1-q}\right) \frac{(1-q^2)dq}{\sqrt{(e_1^2-q^2)(q^2-e_2^2)}} = \int_{1-e_1^2}^{1-e_2^2} \frac{w \operatorname{arctanh}\sqrt{1-w}\,dw}{\sqrt{(1-w)[(1-e_2^2)-w][w-(1-e_1^2)]}}.$$

If we set

$$(2.114)\ \omega = \operatorname{arctanh}\sqrt{1-w},$$

and

$$(2.115)\ d\xi = \frac{w\,dw}{\sqrt{(1-w)[(1-e_2^2)-w][w-(1-e_1^2)]}},$$

then the last integral in Eq. (2.113) can be integrated by parts, and $\xi$ can be expressed in terms of $F(u,k')$ and $E(u,k')$. With the upper limit of the integral for $\xi$ taken as the variable $w$,

$$(2.116)\ \xi_1 = \begin{cases} \displaystyle\int_{1-e_1^2}^{w} \frac{\tilde{w}\,d\tilde{w}}{\sqrt{(1-\tilde{w})[(1-e_2^2)-\tilde{w}][\tilde{w}-(1-e_1^2)]}} \\ 0 \end{cases}$$

$$= \begin{cases} \dfrac{2F(u,k')}{e_1} - 2e_1 E(u,k'), & \text{when}\quad 1 > (1-e_2^2) \geq w > (1-e_1^2), \\ 0, & \text{when}\quad w = (1-e_1^2), \end{cases} \qquad \textbf{GR 251 (3.132.2)}$$

where

$$(2.117)\ u = \arcsin\sqrt{\frac{w-(1-e_1^2)}{e_1^2-e_2^2}} \quad \text{and}\quad k' = \sqrt{1-\frac{e_2^2}{e_1^2}}.$$

Using Eqs. (2.114) and (2.115) and the procedure for integration by parts outlined by Section 2.3, Eq. (2.113) becomes

$$(2.118) \int_{e_2}^{e_1} \ln\left(\frac{1+q}{1-q}\right) \frac{(1-q^2)dq}{\sqrt{(e_1^2-q^2)(q^2-e_2^2)}}$$

$$= \left[\frac{2}{e_1}\mathbf{K}(k') - 2e_1\mathbf{E}(k')\right] \operatorname{arctanh} e_2 + \int_{1-e_1^2}^{1-e_2^2} \left[\frac{F(u,k')}{e_1} - e_1 E(u,k')\right] \frac{dw}{w\sqrt{1-w}}.$$

From Eq. (2.117),

$$(2.119)\ w = (1-e_1^2) + (e_1^2-e_2^2)\sin^2 u,$$

so that, the last integral in Eq. (2.118) can be rewritten in a more familiar form, namely,

$$(2.120) \int_{1-e_1^2}^{1-e_2^2} \left[ \frac{F(u,k')}{e_1} - e_1 E(u,k') \right] \frac{dw}{w\sqrt{1-w}}$$

$$= \frac{2k'^2}{(1-e_1^2)} \int_0^{\pi/2} \frac{F(u,k')\sin u \cos u \, du}{\left(1+k'^2 \sinh^2\mu \sin^2 u\right)\sqrt{1-k'^2 \sin^2 u}} - \frac{2e_1^2 k'^2}{(1-e_1^2)} \int_0^{\pi/2} \frac{E(u,k')\sin u \cos u \, du}{\left(1+k'^2 \sinh^2\mu \sin^2 u\right)\sqrt{1-k'^2 \sin^2 u}}$$

where $\sinh\mu$ is given by Eq. (2.109) and $(u,k')$ are given by Eq. (2.117). Again, using Eqs. (2.32), (2.34), and (2.110) the integral containing $F(u,k')$ in Eq. (2.120) can be isolated. After some algebra, the integral containing $F(u,k')$ can be simplified and reduced to

$$(2.121) \int_0^{\pi/2} \frac{F(u,k')\sin u \cos u \, du}{\left(1+k'^2 \sinh^2\mu \sin^2 u\right)\sqrt{1-k'^2 \sin^2 u}}$$

$$= \frac{-1}{k'^2 \sinh\mu \cosh\mu} \left[ K(k') \operatorname{arctanh}(k \tanh\mu) - \frac{\pi}{2} F(\varphi,k) \right],$$

where

$$(2.122) \ \varphi = \arcsin(\tanh\mu) \text{ and } k' = \sqrt{1-k^2},$$

and $1 > e_1 > e_2 > 0$ implies that $1 > \tanh\mu > 0$ and $1 > k > 0$.

## 2.6 The Integral $I_6(a,b,c) = \int_0^{\pi/2} \frac{F(u,k')\sin u \cos u \, du}{\left(1-k'^2 \cosh^2 v \sin^2 u\right)\sqrt{1-k'^2 \sin^2 u}}$

Following the procedure in Section 2.4, another solution to $\xi$ in Eq. (2.115) can be defined, namely the solution associated with the lower limit of the integral and given by

$$(2.123) \ -\xi_2 = \begin{cases} 0 \\ \int_w^{1-e_2^2} \frac{\tilde{w} d\tilde{w}}{\sqrt{(1-\tilde{w})[(1-e_2^2)-\tilde{w}][\tilde{w}-(1-e_1^2)]}} \end{cases}$$

$$= \begin{cases} 0 & , \text{ when } w = (1-e_2^2), \\ \dfrac{2}{e_1}\left[F(u,k') - e_2^2 \Pi(u,k'^2,k')\right], & \text{ when } 1 > (1-e_2^2) > w \geq (1-e_1^2), \end{cases} \quad \textbf{GR 251 (3.132.3)}$$

where $\Pi(u,k'^2,k')$ is an incomplete elliptic integral of the third kind in the notation of Gradshteyn and Ryzhik[23] and

$$(2.124) \quad u = \arcsin\sqrt{\frac{(1-e_2^2)-w}{k'^2(1-w)}} \text{ and } k' = \sqrt{1-\frac{e_2^2}{e_1^2}}.$$

The elliptic integral $\Pi(u,k'^2,k')$ is a special case, and can be rewritten as

$$(2.125) \quad \Pi(u,k'^2,k') = \frac{e_1^2}{e_2^2}\left[E(u,k') - \frac{\sqrt{[(1-e_2^2)-w][w-(1-e_1^2)]}}{e_1\sqrt{1-w}}\right], \quad \textbf{BY 11 (111.06)}$$

so that, from Eq. (2.123),

(2.126)

$$-\xi_2 = \begin{cases} 0 & , \text{ when } w = (1-e_2^2) \\ \dfrac{2}{e_1}F(u,k') - 2e_1 E(u,k') + \dfrac{2\sqrt{[(1-e_2^2)-w][w-(1-e_1^2)]}}{\sqrt{1-w}}, & \text{ when } 1 > (1-e_2^2) > w \geq (1-e_1^2), \end{cases}$$

where $(u,k')$ are given by Eq. (2.124).

Using Eqs. (2.113), (2.114), (2.115), and (2.126), and the procedure for integration by parts outlined by Section 2.4,

$$(2.127) \quad \int_{e_2}^{e_1} \ln\left(\frac{1+q}{1-q}\right)\frac{(1-q^2)dq}{\sqrt{(e_1^2-q^2)(q^2-e_2^2)}} = \left[\frac{2}{e_1}K(k') - 2e_1 E(k')\right]\operatorname{arctanh} e_1$$

$$+ \int_{1-e_1^2}^{1-e_2^2}\left[e_1 E(u,k') - \frac{F(u,k')}{e_1}\right]\frac{dw}{w\sqrt{1-w}} - \int_{1-e_1^2}^{1-e_2^2}\frac{\sqrt{[(1-e_2^2)-w][w-(1-e_1^2)]}\,dw}{w(1-w)}.$$

Note that the last integral in Eq. (2.127) has been evaluated before in a slightly different form. With $\tau = \sqrt{1-w}$ and Eq. (2.88), the last integral in Eq. (2.127) simplifies to

(2.128)

$$\int_{1-e_1^2}^{1-e_2^2} \frac{\sqrt{[(1-e_1^2)-w][w-(1-e_2^2)]}\, dw}{w(1-w)} = 2\int_{e_2}^{e_1} \frac{\sqrt{(e_1^2-\tau^2)(\tau^2-e_2^2)}\, d\tau}{\tau(1-\tau^2)} = \pi\left[1-e_1 e_2 - \sqrt{(1-e_1^2)(1-e_2^2)}\right].$$

The integral containing $E(u,k')$ and $F(u,k')$ in Eq. (2.127) can be rewritten by inverting $u(w)$ and solving for $w$ from Eq. (2.124), namely,

(2.129) $w = \dfrac{(1-e_2^2) - k'^2 \sin^2 u}{1 - k'^2 \sin^2 u}$,

so that,

(2.130) $\displaystyle\int_{1-e_1^2}^{1-e_2^2} \left[e_1 E(u,k') - \frac{F(u,k')}{e_1}\right] \frac{dw}{w\sqrt{1-w}} = \frac{2e_1 e_2 k'^2}{(1-e_2^2)} \int_0^{\pi/2} \frac{E(u,k')\sin u \cos u\, du}{\left(1-k'^2 \cosh^2 v \sin^2 u\right)\sqrt{1-k'^2 \sin^2 u}}$

$\displaystyle\quad - \frac{2e_2 k'^2}{e_1(1-e_2^2)} \int_0^{\pi/2} \frac{F(u,k')\sin u \cos u\, du}{\left(1-k'^2 \cosh^2 v \sin^2 u\right)\sqrt{1-k'^2 \sin^2 u}}.$

The first integral on the right side of Eq. (2.130) has been evaluated before in Section 2.3, which isolates the second integral as the only integral in Eq. (2.127) to be evaluated. Using Eqs. (2.32), (2.34), (2.94), (2.127), (2.128), and (2.129),

(2.131) $\displaystyle\int_0^{\pi/2} \frac{F(u,k')\sin u \cos u\, du}{\left(1-k'^2 \cosh^2 v \sin^2 u\right)\sqrt{1-k'^2 \sin^2 u}}$

$\displaystyle\qquad = \frac{1}{k'^2 \sinh v \cosh v}\left[\boldsymbol{K}(k')\operatorname{arctanh}\left(\frac{\tanh v}{k}\right) - \frac{\pi}{2} F(\varphi,k)\right],$

where

(2.132) $\varphi = \arcsin\left(\dfrac{\tanh v}{k}\right)$ and $k' = \sqrt{1-k^2}$,

and $1 > e_1 > e_2 > 0$ implies that $1 > k > \tanh v > 0$.

Finally, while it is possible to continue to modify integrals related to the surface area of ellipsoids where $a > b > c > 0$, the six integrals derived and evaluated in this section demonstrate the techniques and the evaluated integrals are useful in their own right.

# 3    Ellipsoids (a,b,c) where c>b>a>0

The surface area of an ellipsoid $(c > b > a > 0)$ projected onto the x-y plane can be written as

$$(3.1) \quad S(c,b,a) = 8 \int_0^a \int_0^{h(x)} \sqrt{1 + \left(\frac{\partial f}{\partial x}\right)^2 + \left(\frac{\partial f}{\partial y}\right)^2} \, dy \, dx ,$$

where

$$(3.2) \quad h(x) = b\sqrt{1 - \frac{x^2}{a^2}} \quad \text{and} \quad f(x,y) = c\sqrt{1 - \frac{x^2}{a^2} - \frac{y^2}{b^2}} ,$$

and the eight-fold symmetry of the ellipsoid has been used to reduce the range of integration. Figs. 3 and 2 illustrate an ellipsoid where $c > b > a > 0$ with its surface area projected onto the x-y plane. Using Eqs. (1.15), (1.18), and (3.2), Eq. (3.1) can be rewritten as

$$(3.3) \quad S(c,b,a) = 8 \int_0^a \int_0^{h(x)} \frac{\sqrt{1 + \frac{\overline{f_1}^2 x^2}{a^2} + \frac{\overline{f_2}^2 y^2}{b^2}}}{\sqrt{1 - \frac{x^2}{a^2} - \frac{y^2}{b^2}}} \, dy \, dx ,$$

where $c > b > a > 0$ implies that $\infty > \overline{f_1} > \overline{f_2} > 0$. Equation (3.3) can be simplified by the same procedure used in Section 2, namely, the consecutive application of scaled elliptic variables, polar coordinates, and change of variables, which results in a familiar form,

$$(3.4) \quad S(c,b,a) = 8ab \int_0^{\pi/2} \int_0^{\pi/2} \sin\theta \sqrt{1 + \overline{q}^2 \sin^2\theta} \, d\theta \, d\phi ,$$

where

$$(3.5) \quad \overline{q}^2 = \overline{f_1}^2 \cos^2\phi + \overline{f_2}^2 \sin^2\phi .$$

Note that $\infty > \overline{f_1} > \overline{f_2} > 0$ and $\frac{\pi}{2} \geq \phi \geq 0$ imply that $\infty > \overline{q}^2 > 0$. Equation (3.4) can be equally written as

$$(3.6) \quad S(c,b,a) = 8ab \int_0^{\pi/2} \int_0^{\pi/2} \sin\theta \sqrt{1 + \overline{f_1}^2 \sin^2\theta} \sqrt{1 - \overline{w}^2 \sin^2\phi} \, d\theta \, d\phi ,$$

where

(3.7) $\quad \bar{w}^2 = \dfrac{(\bar{f}_1^2 - \bar{f}_2^2)\sin^2\theta}{1 + \bar{f}_1^2 \sin^2\theta}.$

Note that $\infty > \bar{f}_1 > \bar{f}_2 > 0$ and $\dfrac{\pi}{2} \geq \theta \geq 0$ imply that $1 > \bar{w}^2 > 0$. Again the limits of $\theta$ and $\phi$ do not depend upon one another, so that Eqs. (3.4) and (3.6) can be integrated in either order. For reference, recall the definition of $(\bar{\varphi}, \bar{k})$ from Eq. (1.20), namely,

(3.8) $\quad \bar{\varphi} = \arctan \bar{f}_1 \text{ and } \bar{k} = \sqrt{1 - \dfrac{\bar{f}_2^2}{\bar{f}_1^2}}.$

### 3.1 The integral $I_1(c,b,a) = \displaystyle\int_0^{\bar{\alpha}} \dfrac{u E(u)\,du}{\left(\bar{k}^2 - u^2\right)^2 \sqrt{\bar{\alpha}^2 - u^2}}$

Equation (3.6) can be easily integrated with respect to $\phi$, which, by definition, reduces to a complete elliptic integral of the second kind, $E(\bar{w})$, so that,

(3.9) $\quad S(c,b,a) = 8ab \displaystyle\int_0^{\pi/2} \sin\theta \sqrt{1 + \bar{f}_1^2 \sin^2\theta}\, E(\bar{w})\, d\theta,$          **GR 181 (2.583.1)**

where $\bar{w}$ is given by Eq. (3.7), assuming the positive square root. For convenience, set $\bar{w} = u$ in Eq. (3.7) and rewrite Eq. (3.9) in terms of $u$, so that, after some algebra,

(3.10) $\quad S(c,b,a) = \dfrac{8ab\bar{\alpha}\bar{k}^2}{\bar{f}_1^2} \displaystyle\int_0^{\bar{\alpha}} \dfrac{u E(u)\,du}{\left(\bar{k}^2 - u^2\right)^2 \sqrt{\bar{\alpha}^2 - u^2}},$

where

(3.11) $\quad \bar{\alpha} = \sqrt{\dfrac{\bar{f}_1^2 - \bar{f}_2^2}{1 + \bar{f}_1^2}}$ and $\bar{k} = \sqrt{1 - \dfrac{\bar{f}_2^2}{\bar{f}_1^2}}.$

Note that $\infty > \bar{f}_1 > \bar{f}_2 > 0$ implies that $1 > \bar{\alpha} > 0$ and $1 > \bar{k} > 0$. Since the surface area of an ellipsoid where $c > b > a > 0$, is known and given by Eq. (1.19), the integral in Eq. (3.10) can be evaluated and is given by

(3.12) $\quad \displaystyle\int_0^{\bar{\alpha}} \dfrac{u E(u)\,du}{\left(\bar{k}^2 - u^2\right)^2 \sqrt{\bar{\alpha}^2 - u^2}} = \dfrac{\pi \bar{f}_1^2}{4\bar{\alpha}\bar{k}^2}\left[\sqrt{\dfrac{1 + \bar{f}_2^2}{1 + \bar{f}_1^2}} + \dfrac{F(\bar{\varphi}, \bar{k})}{\bar{f}_1} + \bar{f}_1 E(\bar{\varphi}, \bar{k})\right],$

where $(\bar{\varphi}, \bar{k})$ are given by Eq. (3.8), and

(3.13) $\dfrac{b}{a} = \sqrt{\dfrac{1+\overline{f}_1^{\,2}}{1+\overline{f}_2^{\,2}}}$

from Eq. (1.21) allows $b$ to be exchanged for $a$ in Eq. (3.10). In terms of $(\overline{\alpha},\overline{k})$, the integral in Eq. (3.12) becomes

(3.14) $\displaystyle\int_0^{\overline{\alpha}} \dfrac{u E(u)\,du}{\left(\overline{k}^2-u^2\right)^2 \sqrt{\overline{\alpha}^2-u^2}} = \dfrac{\pi}{4}\left[\dfrac{\overline{\alpha}\sqrt{1-\overline{\alpha}^2}}{\overline{k}^2\left(\overline{k}^2-\overline{\alpha}^2\right)} + \dfrac{F(\overline{\varphi},\overline{k})}{\overline{k}^2\sqrt{\overline{k}^2-\overline{\alpha}^2}} + \dfrac{\overline{\alpha}^2 E(\overline{\varphi},\overline{k})}{\overline{k}^2\left(\overline{k}^2-\overline{\alpha}^2\right)^{\tfrac{3}{2}}}\right],$

where $\overline{\varphi}$ can be rewritten as $\overline{\varphi} = \arcsin\left(\dfrac{\overline{\alpha}}{\overline{k}}\right)$, and $\infty > \overline{f}_1 > \overline{f}_2 > 0$ together with Eq. (3.11) imply that $1 > \overline{k} > \overline{\alpha} > 0$. In terms of $D(\overline{\varphi},\overline{k}) = \dfrac{1}{\overline{k}^2}\left[F(\overline{\varphi},\overline{k}) - E(\overline{\varphi},\overline{k})\right]$, the integral in Eq. (3.14) becomes

(3.15) $\displaystyle\int_0^{\overline{\alpha}} \dfrac{u E(u)\,du}{\left(\overline{k}^2-u^2\right)^2 \sqrt{\overline{\alpha}^2-u^2}} = \dfrac{\pi}{4}\left\{\dfrac{\overline{\alpha}\sqrt{1-\overline{\alpha}^2}}{\overline{k}^2\left(\overline{k}^2-\overline{\alpha}^2\right)} + \dfrac{\left[F(\overline{\varphi},\overline{k}) - \overline{\alpha}^2 D(\overline{\varphi},\overline{k})\right]}{\left(\overline{k}^2-\overline{\alpha}^2\right)^{\tfrac{3}{2}}}\right\},$

where

(3.16) $\overline{\varphi} = \arcsin\left(\dfrac{\overline{\alpha}}{\overline{k}}\right)$ and $1 > \overline{k} > \overline{\alpha} > 0$.

The closest relative to the integral in Eq. (3.16) found in the literature is repeated below for comparison,

(3.17) $\displaystyle\int_0^{\overline{\alpha}} \dfrac{u E\!\left(\tfrac{u}{\overline{\alpha}}\right)du}{\left(\overline{k}^2-u^2\right)\sqrt{\overline{\alpha}^2-u^2}} = \dfrac{\pi\overline{\alpha}}{2\overline{k}\sqrt{\overline{k}^2-\overline{\alpha}^2}}\left[\boldsymbol{K}\!\left(\dfrac{\overline{\alpha}}{\overline{k}}\right) - \boldsymbol{D}\!\left(\dfrac{\overline{\alpha}}{\overline{k}}\right)\right],$      **PR3 185 (2.16.5.4)**

where $\boldsymbol{D}\!\left(\dfrac{\overline{\alpha}}{\overline{k}}\right) = \dfrac{1}{\overline{k}^2}\left[\boldsymbol{K}\!\left(\dfrac{\overline{\alpha}}{\overline{k}}\right) - \boldsymbol{E}\!\left(\dfrac{\overline{\alpha}}{\overline{k}}\right)\right].$

## 3.2 Integrals that Confirm $S(c,b,a) = 2\pi a^2 \left\{ 1 + \sqrt{\dfrac{1+\bar{f}_1^2}{1+\bar{f}_2^2}} \left[ \dfrac{F(\bar{\varphi},\bar{k})}{\bar{f}_1} + \bar{f}_1 E(\bar{\varphi},\bar{k}) \right] \right\}$

Equation (3.4) can be integrated with respect to $\theta$ in a straightforward manner, and reduces to the expression,

$$(3.18) \quad S(c,b,a) = 8ab \int_0^{\pi/2} \left[ -\dfrac{\cos\theta \sqrt{1+\bar{q}^2 \sin^2\theta}}{2} - \left( \dfrac{1+\bar{q}^2}{2\bar{q}} \right) \arcsin\left( \dfrac{\bar{q}\cos\theta}{\sqrt{1+\bar{q}^2}} \right) \right]_0^{\pi/2} d\phi,$$

or equivalently,

$$(3.19) \quad S(c,b,a) = 2\pi ab + 8ab \int_0^{\pi/2} \left( \dfrac{1+\bar{q}^2}{2\bar{q}} \right) \arcsin\left( \dfrac{\bar{q}}{\sqrt{1+\bar{q}^2}} \right) d\phi, \qquad \textbf{GR 200 (2.598), 182 (2.583.2)}$$

where

$$(3.20) \quad \bar{q}^2 = \bar{f}_1^2 \cos^2\phi + \bar{f}_2^2 \sin^2\phi.$$

Note that $\infty > \bar{f}_1 > \bar{f}_2 > 0$ and $\dfrac{\pi}{2} \geq \phi \geq 0$ imply that $\infty > \bar{q}^2 > 0$. With $\arctan\bar{q} = \arcsin\left( \dfrac{\bar{q}}{\sqrt{1+\bar{q}^2}} \right)$ **[GR 55 (1.624.7)]** and $q(\phi)$ inverted using Eq. (3.20), Eq. (3.19) can be rewritten as a function of $\bar{q}$ only, namely,

$$(3.21) \quad S(c,b,a) = 2\pi ab + 4ab \int_{\bar{f}_2}^{\bar{f}_1} \dfrac{\arctan\bar{q}\, d\bar{q}}{\sqrt{(\bar{f}_1^2 - \bar{q}^2)(\bar{q}^2 - \bar{f}_2^2)}} + 4ab \int_{\bar{f}_2}^{\bar{f}_1} \dfrac{\bar{q}^2 \arctan\bar{q}\, d\bar{q}}{\sqrt{(\bar{f}_1^2 - \bar{q}^2)(\bar{q}^2 - \bar{f}_2^2)}}.$$

Each of the integrals in Eq. (3.21) can be evaluated by variable substitution. Let

$$(3.22) \quad \bar{q} = \tan\bar{\varphi}\sqrt{1 - \bar{k}^2 \sin^2\chi},$$

so that $\chi \in [0, \dfrac{\pi}{2}]$ implies that

$$(3.23) \quad \bar{q} = \left\{ \begin{array}{c} \tan\bar{\varphi} \\ \tan\bar{\varphi}\sqrt{1-\bar{k}^2} \end{array} \right\}, \text{ when } \chi = \left\{ \begin{array}{c} 0 \\ \dfrac{\pi}{2} \end{array} \right\}.$$

Choose $\tan\bar{\varphi} = \bar{f}_1$, so that $\bar{f}_2 = \bar{f}_1\sqrt{1-\bar{k}^2}$ when

(3.24) $\bar{k} = \sqrt{1 - \dfrac{\bar{f}_2^2}{\bar{f}_1^2}}$ .

With these substitutions, the first integral in Eq. (3.21) becomes

(3.25) $\displaystyle\int_{\bar{f}_2}^{\bar{f}_1} \dfrac{\arctan \bar{q}\, d\bar{q}}{\sqrt{(\bar{f}_1^2 - \bar{q}^2)(\bar{q}^2 - \bar{f}_2^2)}} = \dfrac{1}{\tan \bar{\varphi}} \int_0^{\pi/2} \dfrac{\arctan\left(\tan\bar{\varphi}\sqrt{1 - \bar{k}^2 \sin^2\chi}\right) d\chi}{\sqrt{1 - \bar{k}^2 \sin^2\chi}}$ ,

which is proportional to the sum of two integrals given by Gradshteyn and Ryzhik[23], so that,

(3.26) $\displaystyle\int_{\bar{f}_2}^{\bar{f}_1} \dfrac{\arctan \bar{q}\, d\bar{q}}{\sqrt{(\bar{f}_1^2 - \bar{q}^2)(\bar{q}^2 - \bar{f}_2^2)}} = \left(\dfrac{\pi}{2}\right) \dfrac{F(\bar{\varphi}, \bar{k})}{\tan \bar{\varphi}}$ ,   **GR 603 (4.577.1, 2)**

where

(3.27) $\bar{\varphi} = \arctan \bar{f}_1$ and $\bar{k} = \sqrt{1 - \dfrac{\bar{f}_2^2}{\bar{f}_1^2}}$ .

The second integral in Eq. (3.21) can be evaluated by variable substitution in a similar fashion. In this case, let

(3.28) $\bar{q} = \tan \bar{\varphi} \sqrt{1 - \bar{k}^2 x^2}$ ,

so that $x \in [0,1]$ implies that

(3.29) $\bar{q} = \begin{cases} \tan \bar{\varphi} \\ \tan \bar{\varphi} \sqrt{1 - \bar{k}^2} \end{cases}$, when $x = \begin{cases} 0 \\ 1 \end{cases}$.

Again choose $\tan \bar{\varphi} = \bar{f}_1$, so that $\bar{f}_2 = \bar{f}_1 \sqrt{1 - \bar{k}^2}$ when

(3.30) $\bar{k} = \sqrt{1 - \dfrac{\bar{f}_2^2}{\bar{f}_1^2}}$ .

With these substitutions, the second integral in Eq. (3.21) becomes

(3.31) $\displaystyle\int_{\bar{f}_2}^{\bar{f}_1} \dfrac{\bar{q}^2 \arctan \bar{q}\, d\bar{q}}{\sqrt{(\bar{f}_1^2 - \bar{q}^2)(\bar{q}^2 - \bar{f}_2^2)}} = \bar{k} \tan \bar{\varphi} \int_0^1 \left(\dfrac{1 - x^2}{\bar{k}^{-2} - x^2}\right)^{-\tfrac{1}{2}} \arctan\left(\tan \bar{\varphi} \sqrt{1 - \bar{k}^2 x^2}\right) dx$ ,

which is proportional to an integral given by Prudnikov, Brychkov, and Marichev[27, 28, 29], namely,

$$(3.32) \quad \int_{\bar{f}_2}^{\bar{f}_1} \frac{\bar{q}^2 \arctan \bar{q} \, d\bar{q}}{\sqrt{(\bar{f}_1^2 - \bar{q}^2)(\bar{q}^2 - \bar{f}_2^2)}} = \left(\frac{\pi}{2}\right) E(\bar{\varphi}, \bar{k}) \tan \bar{\varphi} - \left(\frac{\pi}{2}\right)\left(1 - \sqrt{1 - \bar{k}^2 \sin^2 \bar{\varphi}}\right), \quad \textbf{PR1 561 (2.7.8.4)}$$

where

$$(3.33) \quad \bar{\varphi} = \arctan \bar{f}_1 \text{ and } \bar{k} = \sqrt{1 - \frac{\bar{f}_2^2}{\bar{f}_1^2}}.$$

Using Eqs. (3.21), (3.26), and (3.32), the surface area becomes

$$(3.34) \quad S(c,b,a) = 2\pi ab + 2\pi ab \left[\frac{F(\bar{\varphi}, \bar{k})}{\tan \bar{\varphi}} + E(\bar{\varphi}, \bar{k}) \tan \bar{\varphi} - \left(1 - \sqrt{1 - \bar{k}^2 \sin^2 \bar{\varphi}}\right)\right],$$

or, in terms of $(\bar{f}_1, \bar{f}_2)$,

$$(3.35) \quad S(c,b,a) = 2\pi a^2 \left\{1 + \sqrt{\frac{1 + \bar{f}_1^2}{1 + \bar{f}_2^2}} \left[\frac{F(\bar{\varphi}, \bar{k})}{\bar{f}_1} + \bar{f}_1 E(\bar{\varphi}, \bar{k})\right]\right\},$$

where

$$(3.36) \quad \bar{\varphi} = \arctan \bar{f}_1 \text{ and } \bar{k} = \sqrt{1 - \frac{\bar{f}_2^2}{\bar{f}_1^2}}.$$

Eq. (3.35) is identical to the expression found in Eq. (1.19). Furthermore, note that Eq. (3.32) could be evaluated from **GR 603 (4.577.1 & .2)** by adding and subtracting with $\bar{k}^2$ multiplying the first integral.

Finally, while Eqs. (3.35) and (3.36) serve as an independent confirmation of the elliptic integral form of the surface area of an ellipsoid where $c > b > a > 0$, the correctness of the expression is not in question.

## 4  Analytic Extensions

In the results derived in Sections 2 and 3, the arguments of functions were assumed to be real numbers, and the limits, intervals, and integration variables of integrals were considered to be real numbers. By analytic extension we mean the process where complementary functions are obtained by allowing their arguments to become complex, and similarly, where complementary integrals are obtained by allowing their limits, intervals, and integration variables to become complex. An illustration of analytic extension, as used in this study, is given for the surface area of an oblate spheroid, where the surface area of the oblate spheroid is analytically extended into the surface area of a prolate spheroid. Discussions of analytic extensions for the integrals $I_N(a,b,c)$, where $N = 1, 2, ..., 6$, and $I_M(c,b,a)$, where $M = 1, 2, 3$, are summarized briefly. A table of results is provided at the end of the section.

To facilitate the study of analytic extensions, a few general equations are referenced, namely,

(4.1) $\quad \begin{Bmatrix} \operatorname{arcsinh}(iz) \\ \operatorname{arccosh} z \\ \operatorname{arctanh}(iz) \\ \operatorname{arccoth}(iz) \end{Bmatrix} = \begin{Bmatrix} i \arcsin z \\ i \arccos z \\ i \arctan z \\ \dfrac{\operatorname{arccot} z}{i} \end{Bmatrix}$, **GR 54 (1.622.1, 2, 3, 4)**

(4.2) $\quad \begin{Bmatrix} \sinh(iz) \\ \cosh(iz) \\ \tanh(iz) \\ \coth(iz) \end{Bmatrix} = \begin{Bmatrix} i \sin z \\ \cos z \\ i \tan z \\ \dfrac{\cot z}{i} \end{Bmatrix}$, **GR 28 (1.311.1, 3, 5, 7)**

(4.3) $\quad \begin{Bmatrix} F(\varphi, ik) \\ E(\varphi, ik) \end{Bmatrix} = \begin{Bmatrix} k_1' F(\beta, k_1) \\ \dfrac{1}{k_1'} \left[ E(\beta, k_1) - \dfrac{k_1^2 \sin \beta \cos \beta}{\sqrt{1 - k_1^2 \sin^2 \beta}} \right] \end{Bmatrix}$, **BY 38 (160.02)**

where

(4.4) $\quad k_1 = \dfrac{k}{\sqrt{1+k^2}}$, $k_1' = \sqrt{1 - k_1^2}$, and $\beta = \arcsin\left( \dfrac{\sqrt{1+k^2} \sin \varphi}{\sqrt{1 + k^2 \sin^2 \varphi}} \right)$, **BY 38 (160.02)**

(4.5) $\quad \begin{Bmatrix} F(i\varphi, k) \\ E(i\varphi, k) \end{Bmatrix} = \begin{Bmatrix} i F(\delta, k') \\ i \left[ F(\delta, k') - E(\delta, k') + \tan \delta \sqrt{1 - k'^2 \sin^2 \delta} \right] \end{Bmatrix}$, **BY 38 (161.02)**

and

(4.6) $\quad \sinh \varphi = \tan \delta$. **BY 38 (161.02)**

In what follows, we shall use Eqs. (4.1) to (4.6) repeatedly without further reference. Two simple examples of analytic extension are considered before beginning with the integrals $I_N(a,b,c)$ and $I_M(c,b,a)$.

The surface areas of an oblate and prolate spheroid are well known and can be found in a number of references[30-32]. Using the convention for $S(a,b,c)$ in Section 2, where $a = b = r > c > 0$, the surface area for an oblate spheroid can be written as

(4.7) $\quad S(r,r,c) = 2\pi r^2 + \dfrac{\pi r c^2}{\sqrt{r^2-c^2}} \ln\left(\dfrac{r+\sqrt{r^2-c^2}}{r-\sqrt{r^2-c^2}}\right).$    **ZW 364 (4.18.1) (r>c)**

Since

(4.8) $\quad \ln\left(\dfrac{r+\sqrt{r^2-c^2}}{r-\sqrt{r^2-c^2}}\right) = 2\,\text{arctanh}\sqrt{\dfrac{r^2-c^2}{r^2}},$    **GR 54 (1.622.7)**

$S(r,r,c)$ can be easily extended to the surface area for a prolate spheroid by using Eq. (4.1) and setting $\sqrt{-z} = i\sqrt{z}$ wherever a negative square root appears. With these changes, Eq. (4.7) becomes

(4.9) $\quad S(r,r,c) = 2\pi r^2 + \dfrac{2\pi r c^2}{\sqrt{c^2-r^2}} \arctan\sqrt{\dfrac{c^2-r^2}{r^2}},$

so that with

(4.10) $\quad \arctan\sqrt{\dfrac{c^2-r^2}{r^2}} = \arcsin\sqrt{\dfrac{c^2-r^2}{c^2}},$    **GR 55 (1.624.7)**

$S(r,r,c)$ becomes

(4.11) $\quad S(r,r,c) = 2\pi r^2 + \dfrac{2\pi r c^2}{\sqrt{c^2-r^2}} \arcsin\sqrt{\dfrac{c^2-r^2}{c^2}},$

which is identical to the surface area of a prolate spheroid, namely,

(4.12) $\quad S(c,r,r) = 2\pi r^2 + \dfrac{2\pi r c^2}{\sqrt{c^2-r^2}} \arcsin\sqrt{\dfrac{c^2-r^2}{c^2}}.$    **ZW 364 (4.18.1) (c>r)**

In the case of an ellipsoid where $a > b > c > 0$, the eccentricities $(e_1, e_2)$ can be rewritten as a imaginary functions of $(\bar{f}_1, \bar{f}_2)$ by reversing the order of the magnitude of $(a,b,c)$ such that $c > b > a > 0$, namely,

(4.13) $\quad e_1 = \dfrac{i\bar{e}_1}{\sqrt{1-\bar{e}_1^2}} = i\bar{f}_1\ \text{ and }\ e_2 = \dfrac{i\bar{e}_2}{\sqrt{1-\bar{e}_2^2}} = i\bar{f}_2,$

where $(\bar{e}_1, \bar{e}_2)$ and $(\bar{f}_1, \bar{f}_2)$ are given by Eqs. (1.15) and (1.18), respectively. The ratio of eccentricities $k$, and its complement $k' = \sqrt{1-k^2}$, then become

(4.14) $k = \dfrac{e_2}{e_1} = \left(\dfrac{\overline{e}_2}{\overline{e}_1}\right)\sqrt{\dfrac{1-\overline{e}_1^2}{1-\overline{e}_2^2}} = \dfrac{\overline{f}_2}{\overline{f}_1}$ and $k' = \sqrt{1-\dfrac{e_2^2}{e_1^2}} = \sqrt{\dfrac{\overline{e}_1^2 - \overline{e}_2^2}{\overline{e}_1^2(1-\overline{e}_2^2)}} = \sqrt{1-\dfrac{\overline{f}_2^2}{\overline{f}_1^2}}$.

Note that from Equation (1.20), $\overline{k} = k'$ and $\overline{k}' = k$, where $\infty > \overline{f}_1 > \overline{f}_2 > 0$ implies that $1 > (k, k', \overline{k}, \overline{k}') > 0$, so that the incomplete elliptic integrals found in Sections 2 and 3 will not need to be analytically extended. However, the modulus of complete elliptic integrals and the argument of incomplete elliptic integrals do become purely imaginary when $a > b > c > 0$ inverts to $c > b > a > 0$. These are considered on a case-by-case basis using Eqs. (4.3) through (4.6) above.

As a second example, consider Eq. (2.32), which is repeated below for reference,

(4.15) $\displaystyle\int_{e_2}^{e_1} \ln\left(\dfrac{1+q}{1-q}\right) \dfrac{dq}{\sqrt{e_1^2 - q^2}\sqrt{q^2 - e_2^2}} = \dfrac{\pi}{e_1} F(\varphi, k)$,    **PR1 510 (2.6.13.2)**

where

(4.16) $\varphi = \arcsin e_1$ and $k = \dfrac{e_2}{e_1}$.

Note that at this point it is important to identify and separate square roots in Eq. (4.15) so that the associated branch cuts can be properly considered and analytically extended. In short, $\sqrt{xy} = \sqrt{x}\sqrt{y}$ only when $x$ and $y$ are non-negative.

With $\ln\left(\dfrac{1+q}{1-q}\right) = 2\operatorname{arctanh} q$ [**GR 54 (1.622.7)**], the analytic extension, $e_1 = i\overline{f}_1$ and $e_2 = i\overline{f}_2$, and subsequent transformation, $q = i\overline{q}$ and $\varphi = i\overline{\varphi}$, reduce Eq. (4.15) to

(4.17) $\displaystyle\int_{\overline{f}_2}^{\overline{f}_1} \dfrac{\arctan \overline{q}\, d\overline{q}}{\sqrt{\overline{f}_1^2 - \overline{q}^2}\sqrt{\overline{q}^2 - \overline{f}_2^2}} = \int_{\overline{f}_2}^{\overline{f}_1} \dfrac{\arctan \overline{q}\, d\overline{q}}{\sqrt{(\overline{f}_1^2 - \overline{q}^2)(\overline{q}^2 - \overline{f}_2^2)}} = \dfrac{\pi}{2\overline{f}_1} F(\overline{\varphi}, \overline{k})$,

where $\infty > \overline{f}_1 > \overline{f}_2 > 0$, and

(4.18) $\overline{\varphi} = \arctan \overline{f}_1$ and $\overline{k} = \sqrt{1 - \dfrac{\overline{f}_2^2}{\overline{f}_1^2}}$.

Since Eqs. (4.17) and (4.18) are identical to Eqs. (3.26) and (3.27), the analytic extension of Eq. (4.15) corresponds to the sum of integrals given by Gradshteyn and Ryzhik[23], namely, **GR 603 (4.577.1)** and **GR 603 (4.577.2)**.

## 4.1 Extension of $I_1(a,b,c) = \int_0^{\bar{\alpha}} \dfrac{uE(u)du}{\left(k'^2 + k^2 u^2\right)^2 \sqrt{\bar{\alpha}^2 - u^2}}$

Consider the integral in Eq. (2.21), repeated below for reference,

(4.19) $\displaystyle\int_0^{\alpha} \dfrac{uE(u)du}{\left(k'^2 + k^2 u^2\right)^2 \sqrt{\alpha^2 - u^2}} = \dfrac{\pi}{4}\left[\dfrac{\alpha\sqrt{1-\alpha^2}}{\left(k'^2 + k^2\alpha^2\right)^2} + \dfrac{\alpha^2 E(\lambda,k)}{k'^2\left(k'^2 + k^2\alpha^2\right)^{\frac{3}{2}}} + \dfrac{\left(1-\alpha^2\right)F(\lambda,k)}{\left(k'^2 + k^2\alpha^2\right)^{\frac{3}{2}}}\right],$

where $(\lambda, k, k', \alpha)$ are given by Eq. (2.22), and note that the analytic extension, $e_1 = i\bar{f}_1$ and $e_2 = i\bar{f}_2$, modifies $(\alpha, k, k')$ such that,

(4.20) $\alpha = \dfrac{i\bar{\alpha}}{\sqrt{1-\bar{\alpha}^2}}$, $k = \bar{k}'$, and $k' = \bar{k}$,

where $(\bar{\alpha}, \bar{k})$ are given by Eq. (3.11) and $\bar{k}' = \sqrt{1-\bar{k}^2}$. With $u = \dfrac{ix}{\sqrt{1-x^2}}$, the integral in Eq. (4.19) becomes,

(4.21) $\displaystyle\int_0^{\alpha} \dfrac{uE(u)du}{\left(k'^2 + k^2 u^2\right)^2 \sqrt{\alpha^2 - u^2}} = i\sqrt{1-\bar{\alpha}^2}\int_0^{\bar{\alpha}} \dfrac{xE(x)dx}{\left(\bar{k}^2 - x^2\right)^2 \sqrt{\bar{\alpha}^2 - x^2}},$

so that, with $\lambda = i\bar{\varphi}$, the remaining terms in Eq. (4.19) become,

(4.22) $\dfrac{\pi}{4}\left[\dfrac{\alpha\sqrt{1-\alpha^2}}{\left(k'^2 + k^2\alpha^2\right)^2} + \dfrac{\alpha^2 E(\lambda,k)}{k'^2\left(k'^2 + k^2\alpha^2\right)^{\frac{3}{2}}} + \dfrac{\left(1-\alpha^2\right)F(\lambda,k)}{\left(k'^2 + k^2\alpha^2\right)^{\frac{3}{2}}}\right]$

$= \dfrac{i\pi}{4}\left[\dfrac{\bar{\alpha}\left(1-\bar{\alpha}^2\right)}{\bar{k}^2\left(\bar{k}^2 - \bar{\alpha}^2\right)} + \dfrac{\sqrt{1-\bar{\alpha}^2}\,F(\bar{\varphi},\bar{k})}{\bar{k}^2\sqrt{\bar{k}^2 - \bar{\alpha}^2}} + \dfrac{\bar{\alpha}^2\sqrt{1-\bar{\alpha}^2}\,E(\bar{\varphi},\bar{k})}{\bar{k}^2\left(\bar{k}^2 - \bar{\alpha}^2\right)^{\frac{3}{2}}}\right],$

so that, together, Eqs. (4.21) and (4.22) imply that,

(4.23) $\displaystyle\int_0^{\bar{\alpha}} \dfrac{xE(x)dx}{\left(\bar{k}^2 - x^2\right)^2 \sqrt{\bar{\alpha}^2 - x^2}} = \dfrac{\pi}{4}\left\{\dfrac{\bar{\alpha}\sqrt{1-\bar{\alpha}^2}}{\bar{k}^2\left(\bar{k}^2 - \bar{\alpha}^2\right)} + \dfrac{F(\bar{\varphi},\bar{k})}{\bar{k}^2\sqrt{\bar{k}^2 - \bar{\alpha}^2}} + \dfrac{\bar{\alpha}^2 E(\bar{\varphi},\bar{k})}{\bar{k}^2\left(\bar{k}^2 - \bar{\alpha}^2\right)^{\frac{3}{2}}}\right\},$

where

(4.24) $\bar{\varphi} = \arcsin\left(\dfrac{\bar{\alpha}}{\bar{k}}\right)$ and $1 > \bar{k} > \bar{\alpha} > 0$.

The inequality arises from $\infty > \bar{f}_1 > \bar{f}_2 > 0$ together with Eq. (3.11). Since Eqs. (4.23) and (4.24) are identical to Eqs. (3.14) and (3.16), $I_1(a,b,c)$ is the analytic extension of $I_1(c,b,a)$.

## 4.2 Extension of $I_2(a,b,c) = \displaystyle\int_\alpha^\beta \ln\left(\dfrac{\varepsilon+u}{\varepsilon-u}\right)\dfrac{u^2\,du}{\sqrt{(u^2-\alpha^2)(\beta^2-u^2)}}$

With Eq. (2.34) repeated below for reference,

(4.25) $\displaystyle\int_{e_2}^{e_1} \ln\left(\dfrac{1+q}{1-q}\right)\dfrac{q^2\,dq}{\sqrt{(e_1^2-q^2)(q^2-e_2^2)}} = \pi\left[1 - \sqrt{(1-e_1^2)(1-e_2^2)}\right] + \pi e_1\left[F(\varphi,k) - E(\varphi,k)\right],$

where $(\varphi, k)$ are given by Eq. (2.33), and the procedures for analytic extension described earlier, Eq. (4.25) can be reduced to a familiar form by $e_1 = i\bar{f}_1$ and $e_2 = i\bar{f}_2$, and subsequently $q = i\bar{q}$, namely, Eqs. (3.32) and (3.33), expressed below in terms of $(\bar{f}_1, \bar{f}_2)$,

(4.26) $\displaystyle\int_{\bar{f}_2}^{\bar{f}_1} \dfrac{\bar{q}^2 \arctan\bar{q}\,d\bar{q}}{\sqrt{(\bar{f}_1^2-\bar{q}^2)(\bar{q}^2-\bar{f}_2^2)}} = \left(\dfrac{\pi}{2}\right)\bar{f}_1 E(\bar{\varphi},\bar{k}) - \left(\dfrac{\pi}{2}\right)\left(1 - \sqrt{\dfrac{1+\bar{f}_2^2}{1+\bar{f}_1^2}}\right),$

where

(4.27) $\bar{\varphi} = \arctan\bar{f}_1$ and $\bar{k} = \sqrt{1 - \dfrac{\bar{f}_2^2}{\bar{f}_1^2}}$.

With the variable substitution given by Eq. (3.28), namely, $\bar{q} = \tan\bar{\varphi}\sqrt{1-\bar{k}^2 x^2}$, where $\tan\bar{\varphi} = \bar{f}_1$ and $\bar{f}_2 = \bar{f}_1\sqrt{1-\bar{k}^2}$, the integral in Eq. (4.26) can be rewritten as

(4.28) $\displaystyle\int_{\bar{f}_2}^{\bar{f}_1} \dfrac{\bar{q}^2 \arctan\bar{q}\,d\bar{q}}{\sqrt{(\bar{f}_1^2-\bar{q}^2)(\bar{q}^2-\bar{f}_2^2)}} = \bar{k}\tan\bar{\varphi}\int_0^1 \left(\dfrac{1-x^2}{\bar{k}^{-2}-x^2}\right)^{-\frac{1}{2}} \arctan\left(\tan\bar{\varphi}\sqrt{1-\bar{k}^2 x^2}\right)dx,$

which is proportional to an integral given by Prudnikov, Brychkov, and Marichev[27-29], namely, **PR1 561 (2.7.8.4).** This implies that the integral in Eq. (4.26) could be equally written as Eq. (3.32).

With similar techniques, $I_2(a,b,c)$, given by Eqs. (2.74) and (2.75), can be reduced to Eq. (4.28) by the variable substitution $\alpha = i\varepsilon\bar{f}_2$ and $\beta = i\varepsilon\bar{f}_1$, and subsequently $u = i\varepsilon\bar{q}$, so that $I_2(a,b,c)$ is the analytic extension of an integral which is similarly proportional to an integral given by Prudnikov, Brychkov, and Marichev[27-29], namely, **PR1 561 (2.7.8.4).**

### 4.3  Extension of $I_3(a,b,c) = \int_0^{\pi/2} \dfrac{E(u,k')\sin u \cos u \, du}{\left(1-k'^2 \cosh^2 v \sin^2 u\right)\sqrt{1-k'^2 \sin^2 u}}$

The analytic extension of $I_3(a,b,c)$ is a bit more complicated than others before, in the sense that there does not appear to be a counterpart in the literature. The analytic extension then, results in yet another integral expression, here referred to simply as $I_2(c,b,a)$. With Eq. (2.94) repeated below for reference,

$$(4.29) \quad \int_0^{\pi/2} \frac{E(u,k')\sin u \cos u \, du}{\left(1-k'^2 \cosh^2 v \sin^2 u\right)\sqrt{1-k'^2 \sin^2 u}}$$

$$= \frac{1}{k'^2 \sinh v \cosh v}\left\{E(k')\operatorname{arctanh}\left(\frac{\tanh v}{k}\right) - \frac{\pi \tanh v}{2} - \frac{\pi}{2}\left[F(\varphi,k) - E(\varphi,k)\right]\right\},$$

where $(\varphi, k')$ are given by Eq. (2.95) and $1 > k > \tanh v > 0$, the analytic extension, $e_1 = i\overline{f_1}$ and $e_2 = i\overline{f_2}$, implies that $v$ is purely imaginary from Eq. (2.93), so that with $v = i\psi$, Eq. (4.29) reduces to

$$(4.30) \quad \int_0^{\pi/2} \frac{E(u,\overline{k})\sin u \cos u \, du}{\left(1-\overline{k}^2 \cos^2 \psi \sin^2 u\right)\sqrt{1-\overline{k}^2 \sin^2 u}}$$

$$= \frac{1}{\overline{k}^2 \sin\psi \cos\psi}\left[E(\overline{k})\arctan\left(\frac{\tan\psi}{\overline{k}'}\right) - \frac{\pi}{2}E(\beta,\overline{k}) + \frac{\pi}{2}\frac{\tan\psi}{\sqrt{1-\overline{k}^2 \cos^2\psi}}\left(1-\sqrt{1-\overline{k}^2 \cos^2\psi}\right)\right],$$

where

$$(4.31) \quad \beta = \arctan\left(\frac{\tan\psi}{\overline{k}'}\right) \text{ and } \overline{k}' = \sqrt{1-\overline{k}^2}, \text{ and } 1 > \overline{k} > 0 \text{ and } \frac{\pi}{2} > \psi > 0.$$

While Eqs. (4.30) and (4.31) do not appear to have a counterpart in the literature per se, it is possible to convert Eq. (4.30) into a $\sin^2\psi$ version by using the trigonometric formulas,

$$(4.32) \quad \left\{\begin{array}{l}\sin\left(x\pm\dfrac{\pi}{2}\right)\\[4pt]\cos\left(x\pm\dfrac{\pi}{2}\right)\end{array}\right\} = \left\{\begin{array}{l}\pm\cos x\\ \mp\sin x\end{array}\right\}, \qquad\qquad\qquad \textbf{GR 29 (1.313.1, 5)}$$

and the addition formula for incomplete elliptic integrals,

(4.33) $\begin{Bmatrix} F(\vartheta,k)+F(\delta,k) \\ E(\vartheta,k)+E(\delta,k) \end{Bmatrix} = \begin{Bmatrix} K(k) \\ E(k)+k^2 \sin\vartheta\sin\delta \end{Bmatrix}$, **BY 13 (117.01)**

where

(4.34) $\cot\delta = k'\tan\vartheta$.

The algebra is a bit tedious and care must be taken to keep functions within their principal values, but the end result is the one expected, namely, that Eqs. (4.30) and (4.31) ultimately reduce to an integral given by Gradshteyn and Ryzhik[23], namely, **GR 626 (6.123)**. Without reproducing that algebra, the results may be summarized; that is, the analytic extension of $I_3(a,b,c)$ is $I_2(c,b,a)$, given by Eqs. (4.30) and (4.31), and $I_2(c,b,a)$ is the $\cos^2\psi$ version of **GR 626 (6.123)**.

**4.4 Extension of** $I_4(a,b,c) = \int_0^{\pi/2} \dfrac{E(u,k')\sin u \cos u \, du}{\left(1+k'^2\sinh^2\mu\sin^2 u\right)\sqrt{1-k'^2\sin^2 u}}$

Equation (2.110) is repeated below for reference,

(4.35) $\int_0^{\pi/2} \dfrac{E(u,k')\sin u \cos u \, du}{\left(1+k'^2\sinh^2\mu\sin^2 u\right)\sqrt{1-k'^2\sin^2 u}}$

$= \dfrac{-1}{k'^2 \sinh\mu\cosh\mu} \left\{ E(k')\operatorname{arctanh}(k\tanh\mu) - \dfrac{\pi}{2}\left[F(\varphi,k)-E(\varphi,k)+\tanh\mu\sqrt{1+k'^2\sinh^2\mu}\right] \right.$

$\left. -\dfrac{\pi}{2}\coth\mu\left(1-\sqrt{1+k'^2\sinh^2\mu}\right) \right\}$,

where $(\varphi,k')$ are given by Eq. (2.111), and $(\mu,k)$ satisfy $1>\tanh\mu>0$ and $1>k>0$. With the analytic extension procedures used earlier for $e_1 = i\overline{f_1}$ and $e_2 = i\overline{f_2}$, and $\mu = i\xi$, Eq. (4.35) can be reduced to an integral given by Gradshteyn and Ryzhik[23], namely,

(4.36) $\int_0^{\pi/2} \dfrac{E(u,\overline{k})\sin u \cos u \, du}{\left(1-\overline{k}^2\sin^2\xi\sin^2 u\right)\sqrt{1-\overline{k}^2\sin^2 u}}$ **GR 626 (6.123)**

$= \dfrac{-1}{\overline{k}^2 \sin\xi\cos\xi}\left[E(\overline{k})\arctan(\overline{k}'\tan\xi) - \dfrac{\pi}{2}E(\xi,\overline{k}) + \dfrac{\pi}{2}\cot\xi\left(1-\sqrt{1-\overline{k}^2\sin^2\xi}\right)\right]$,

where $\bar{k}' = \sqrt{1-\bar{k}^2}$, and $1 > \bar{k} > 0$ and $\frac{\pi}{2} > \xi > 0$. Consequently, the analytic extension of $I_4(a,b,c)$ is the integral **GR 626 (6.123)**.

### 4.5 Extension of $I_5(a,b,c) = \int_0^{\pi/2} \frac{F(u,k')\sin u \cos u\, du}{\left(1+k'^2 \sinh^2\mu \sin^2 u\right)\sqrt{1-k'^2 \sin^2 u}}$

Equation (2.121) is repeated below for reference,

$$(4.37) \quad \int_0^{\pi/2} \frac{F(u,k')\sin u \cos u\, du}{\left(1+k'^2 \sinh^2\mu \sin^2 u\right)\sqrt{1-k'^2 \sin^2 u}}$$

$$= \frac{-1}{k'^2 \sinh\mu \cosh\mu}\left[K(k')\operatorname{arctanh}(k \tanh\mu) - \frac{\pi}{2}F(\varphi,k)\right],$$

where $(\varphi,k')$ are given by Eq. (2.122), and $(\mu,k)$ satisfy $1 > \tanh\mu > 0$ and $1 > k > 0$. With the analytic extension procedures used earlier for $e_1 = i\bar{f_1}$ and $e_2 = i\bar{f_2}$, and $\mu = i\xi$, Eq. (4.37) can be reduced to an integral given by Gradshteyn and Ryzhik[23], namely,

$$(4.38) \quad \int_0^{\pi/2} \frac{F(u,\bar{k})\sin u \cos u\, du}{\left(1-\bar{k}^2 \sin^2\xi \sin^2 u\right)\sqrt{1-\bar{k}^2 \sin^2 u}} \qquad \textbf{GR 625 (6.113.2)}$$

$$= \frac{-1}{\bar{k}^2 \sin\xi \cos\xi}\left[K(\bar{k})\arctan\left(\bar{k}'\tan\xi\right) - \frac{\pi}{2}F(\xi,\bar{k})\right],$$

where $\bar{k}' = \sqrt{1-\bar{k}^2}$, and $1 > \bar{k} > 0$ and $\frac{\pi}{2} > \xi > 0$. Consequently, the analytic extension of $I_5(a,b,c)$ is the integral **GR 625 (6.113.2)**.

### 4.6 Extension of $I_6(a,b,c) = \int_0^{\pi/2} \frac{F(u,k')\sin u \cos u\, du}{\left(1-k'^2 \cosh^2 v \sin^2 u\right)\sqrt{1-k'^2 \sin^2 u}}$

The analytic extension of $I_6(a,b,c)$ is similar to $I_3(a,b,c)$ in the sense that it is more complicated and there does not appear to be a counterpart in the literature. The analytic extension then, results in yet another integral expression, here referred to simply as $I_3(c,b,a)$. With Eq. (2.131) repeated below for reference,

(4.39) $$\int_0^{\pi/2} \frac{F(u,k')\sin u \cos u \, du}{(1-k'^2 \cosh^2 v \sin^2 u)\sqrt{1-k'^2 \sin^2 u}}$$

$$= \frac{1}{k'^2 \sinh v \cosh v}\left[ K(k')\operatorname{arctanh}\left(\frac{\tanh v}{k}\right) - \frac{\pi}{2}F(\varphi,k)\right],$$

where $(\varphi, k')$ are given by Eq. (2.132) and $1 > k > \tanh v > 0$, the analytic extension procedures used earlier for $e_1 = i f_1$ and $e_2 = i f_2$, and $v = i\psi$, Eq. (4.39) can be reduced to another integral expression, namely,

(4.40) $$\int_0^{\pi/2} \frac{F(u,\bar{k})\sin u \cos u \, du}{(1-\bar{k}^2 \cos^2\psi \sin^2 u)\sqrt{1-\bar{k}^2 \sin^2 u}} = \frac{1}{\bar{k}^2 \sin\psi \cos\psi}\left[ K(\bar{k})\arctan\left(\frac{\tan\psi}{\bar{k}'}\right) - \frac{\pi}{2}F(\beta,\bar{k})\right],$$

where

(4.41) $\beta = \arctan\left(\dfrac{\tan\psi}{\bar{k}'}\right)$ and $\bar{k}' = \sqrt{1-\bar{k}^2}$, and $1 > \bar{k} > 0$ and $\dfrac{\pi}{2} > \psi > 0$.

Once again, Eqs. (4.40) and (4.41) do not appear to have a counterpart in the literature per se, but, once again, it is possible to convert Eq. (4.40) into a $\sin^2\psi$ version by using Eqs. (4.32) through (4.34). The end result is the one expected, namely that Eqs. (4.40) and (4.41) ultimately reduce to an integral given by Gradshteyn and Ryzhik[23],

(4.42) $$\int_0^{\pi/2} \frac{F(u,\bar{k})\sin u \cos u \, du}{(1-\bar{k}^2 \sin^2\psi \sin^2 u)\sqrt{1-\bar{k}^2 \sin^2 u}} \qquad \textbf{GR 625 (6.113.2)}$$

$$= \frac{-1}{\bar{k}^2 \sin\psi \cos\psi}\left[ K(\bar{k})\arctan\left(\bar{k}' \tan\psi\right) - \frac{\pi}{2}F(\psi,\bar{k})\right],$$

where $\bar{k}' = \sqrt{1-\bar{k}^2}$, and $1 > \bar{k} > 0$ and $\dfrac{\pi}{2} > \psi > 0$. In summary, the analytic extension of $I_6(a,b,c)$ is $I_3(c,b,a)$, given by Eqs. (4.40) and (4.41), and $I_3(c,b,a)$ is the $\cos^2\psi$ version of **GR 625 (6.113.2)**.

## 4.7 Table of Integrals and Analytic Extensions

| Identifier | Integral | Extension | Comments |
|---|---|---|---|
| $I_1(a,b,c)$ | $\displaystyle\int_0^{\alpha} \frac{u E(u)\,du}{\left(k'^2 + k^2 u^2\right)^2 \sqrt{\alpha^2 - u^2}}$ | $\displaystyle\int_0^{\bar{\alpha}} \frac{u E(u)\,du}{\left(\bar{k}^2 - u^2\right)^2 \sqrt{\bar{\alpha}^2 - u^2}}$ | **Neither integral nor extension found in the literature.** |
| $I_2(a,b,c)$ | $\displaystyle\int_{\alpha}^{\beta} \ln\!\left(\frac{\varepsilon+u}{\varepsilon-u}\right) \frac{u^2\,du}{\sqrt{(u^2-\alpha^2)(\beta^2-u^2)}}$ | $\displaystyle\int_0^{1} \left(\frac{1-x^2}{\bar{k}^{-2}-x^2}\right)^{-\frac{1}{2}} \arctan\!\left(\tan\bar{\varphi}\sqrt{1-\bar{k}^2 x^2}\right) dx$ | **Integral not found in the literature. Extension related to PR1 561 (2.7.8.4).** |
| $I_3(a,b,c)$ | $\displaystyle\int_0^{\pi/2} \frac{E(u,k')\sin u \cos u\,du}{\left(1-k'^2\cosh^2\! v\,\sin^2 u\right)\sqrt{1-k'^2\sin^2 u}}$ | $\displaystyle\int_0^{\pi/2} \frac{E(u,\bar{k})\sin u \cos u\,du}{\left(1-\bar{k}^2\cos^2\!\psi\,\sin^2 u\right)\sqrt{1-\bar{k}^2\sin^2 u}}$ | **Neither integral nor extension found in the literature. Extension converted to GR 626 (6.123) by trigonometric and elliptic integral addition formulas.** |
| $I_4(a,b,c)$ | $\displaystyle\int_0^{\pi/2} \frac{E(u,k')\sin u \cos u\,du}{\left(1+k'^2\sinh^2\!\mu\,\sin^2 u\right)\sqrt{1-k'^2\sin^2 u}}$ | $\displaystyle\int_0^{\pi/2} \frac{E(u,\bar{k})\sin u \cos u\,du}{\left(1-\bar{k}^2\sin^2\!\xi\,\sin^2 u\right)\sqrt{1-\bar{k}^2\sin^2 u}}$ | **Integral not found in the literature. Extension given by GR 626 (6.123)** |
| $I_5(a,b,c)$ | $\displaystyle\int_0^{\pi/2} \frac{F(u,k')\sin u \cos u\,du}{\left(1+k'^2\sinh^2\!\mu\,\sin^2 u\right)\sqrt{1-k'^2\sin^2 u}}$ | $\displaystyle\int_0^{\pi/2} \frac{F(u,\bar{k})\sin u \cos u\,du}{\left(1-\bar{k}^2\sin^2\!\xi\,\sin^2 u\right)\sqrt{1-\bar{k}^2\sin^2 u}}$ | **Integral not found in the literature. Extension given by GR 625 (6.113.2)** |
| $I_6(a,b,c)$ | $\displaystyle\int_0^{\pi/2} \frac{F(u,k')\sin u \cos u\,du}{\left(1-k'^2\cosh^2\! v\,\sin^2 u\right)\sqrt{1-k'^2\sin^2 u}}$ | $\displaystyle\int_0^{\pi/2} \frac{F(u,\bar{k})\sin u \cos u\,du}{\left(1-\bar{k}^2\cos^2\!\psi\,\sin^2 u\right)\sqrt{1-\bar{k}^2\sin^2 u}}$ | **Neither integral nor extension found in the literature. Extension converted to GR 625 (6.113.2) by trigonometric and elliptic integral addition formulas.** |

## 4.8 Summary of Integrals and Series

The main contributions from this study of the surface area of arbitrary ellipsoids fall into three categories as indicated below.

### 4.8.1 Integrals Containing $E(u)$

(4.43) $\displaystyle\int_0^\alpha \frac{u E(u)\,du}{\left(k'^2+k^2u^2\right)^2\sqrt{\alpha^2-u^2}} = \frac{\pi}{4}\left[\frac{\alpha\sqrt{1-\alpha^2}}{\left(k'^2+k^2\alpha^2\right)^2} + \frac{\alpha^2 E(\lambda,k)}{k'^2\left(k'^2+k^2\alpha^2\right)^{\frac{3}{2}}} + \frac{\left(1-\alpha^2\right)F(\lambda,k)}{\left(k'^2+k^2\alpha^2\right)^{\frac{3}{2}}}\right],$

where $\lambda = \arcsin\left(\dfrac{\alpha}{\sqrt{k'^2+k^2\alpha^2}}\right)$ and $k' = \sqrt{1-k^2}$, and $1 > k > 0$ and $1 > \alpha > 0$.

(4.44) $\displaystyle\int_0^{\bar{\alpha}} \frac{u E(u)\,du}{\left(\bar{k}^2-u^2\right)^2\sqrt{\bar{\alpha}^2-u^2}} = \frac{\pi}{4}\left[\frac{\bar{\alpha}\sqrt{1-\bar{\alpha}^2}}{\bar{k}^2\left(\bar{k}^2-\bar{\alpha}^2\right)} + \frac{F(\bar{\varphi},\bar{k})}{\bar{k}^2\sqrt{\bar{k}^2-\bar{\alpha}^2}} + \frac{\bar{\alpha}^2 E(\bar{\varphi},\bar{k})}{\bar{k}^2\left(\bar{k}^2-\bar{\alpha}^2\right)^{\frac{3}{2}}}\right],$

where $\bar{\varphi} = \arcsin\left(\dfrac{\bar{\alpha}}{\bar{k}}\right)$ and $1 > \bar{k} > \bar{\alpha} > 0$.

### 4.8.2 Integrals Containing $E(u,\kappa)$ and $F(u,\kappa)$

(4.45) $\displaystyle\int_0^{\pi/2} \frac{E(u,k')\sin u \cos u\,du}{\left(1-k'^2\cosh^2 v \sin^2 u\right)\sqrt{1-k'^2\sin^2 u}}$

$= \dfrac{1}{k'^2\sinh v \cosh v}\left\{E(k')\operatorname{arctanh}\left(\dfrac{\tanh v}{k}\right) - \dfrac{\pi \tanh v}{2} - \dfrac{\pi}{2}\left[F(\varphi,k)-E(\varphi,k)\right]\right\},$

where $\varphi = \arcsin\left(\dfrac{\tanh v}{k}\right)$ and $k' = \sqrt{1-k^2}$, and $1 > k > \tanh v > 0$.

(4.46) $\displaystyle\int_0^{\pi/2} \frac{E(u,\bar{k})\sin u \cos u\,du}{\left(1-\bar{k}^2\cos^2\psi \sin^2 u\right)\sqrt{1-\bar{k}^2\sin^2 u}}$

$= \dfrac{1}{\bar{k}^2\sin\psi\cos\psi}\left[E(\bar{k})\arctan\left(\dfrac{\tan\psi}{\bar{k}'}\right) - \dfrac{\pi}{2}E(\beta,\bar{k}) + \dfrac{\pi}{2}\dfrac{\tan\psi}{\sqrt{1-\bar{k}^2\cos^2\psi}}\left(1-\sqrt{1-\bar{k}^2\cos^2\psi}\right)\right],$

where $\beta = \arctan\left(\dfrac{\tan\psi}{\bar{k}'}\right)$ and $\bar{k}' = \sqrt{1-\bar{k}^2}$, and $1 > \bar{k} > 0$ and $\dfrac{\pi}{2} > \psi > 0$.

(4.47) $\int_0^{\pi/2} \dfrac{E(u,k')\sin u \cos u \, du}{\left(1+k'^2 \sinh^2\mu \sin^2 u\right)\sqrt{1-k'^2 \sin^2 u}}$

$= \dfrac{-1}{k'^2 \sinh\mu \cosh\mu}\left\{ E(k')\operatorname{arctanh}(k\tanh\mu) - \dfrac{\pi}{2}\left[F(\varphi,k) - E(\varphi,k) + \tanh\mu\sqrt{1+k'^2 \sinh^2\mu}\right]\right.$

$\left. -\dfrac{\pi}{2}\coth\mu\left(1-\sqrt{1+k'^2 \sinh^2\mu}\right)\right\}$,

where $\varphi = \arcsin(\tanh\mu)$ and $k' = \sqrt{1-k^2}$, and $1 > k > 0$ and $1 > \tanh\mu > 0$.

(4.48) $\int_0^{\pi/2} \dfrac{F(u,k')\sin u \cos u \, du}{\left(1+k'^2 \sinh^2\mu \sin^2 u\right)\sqrt{1-k'^2 \sin^2 u}}$

$= \dfrac{-1}{k'^2 \sinh\mu \cosh\mu}\left[ K(k')\operatorname{arctanh}(k\tanh\mu) - \dfrac{\pi}{2}F(\varphi,k)\right]$,

where $\varphi = \arcsin(\tanh\mu)$ and $k' = \sqrt{1-k^2}$, and $1 > k > 0$ and $1 > \tanh\mu > 0$.

(4.49) $\int_0^{\pi/2} \dfrac{F(u,k')\sin u \cos u \, du}{\left(1-k'^2 \cosh^2\nu \sin^2 u\right)\sqrt{1-k'^2 \sin^2 u}}$

$= \dfrac{1}{k'^2 \sinh\nu \cosh\nu}\left[ K(k')\operatorname{arctanh}\left(\dfrac{\tanh\nu}{k}\right) - \dfrac{\pi}{2}F(\varphi,k)\right]$,

where $\varphi = \arcsin\left(\dfrac{\tanh\nu}{k}\right)$ and $k' = \sqrt{1-k^2}$, and $1 > k > \tanh\nu > 0$.

(4.50) $\int_0^{\pi/2} \dfrac{F(u,\bar{k})\sin u \cos u \, du}{\left(1-\bar{k}^2 \cos^2\psi \sin^2 u\right)\sqrt{1-\bar{k}^2 \sin^2 u}} = \dfrac{1}{\bar{k}^2 \sin\psi \cos\psi}\left[ K(\bar{k})\arctan\left(\dfrac{\tan\psi}{\bar{k}'}\right) - \dfrac{\pi}{2}F(\beta,\bar{k})\right]$,

where $\beta = \arctan\left(\dfrac{\tan\psi}{\bar{k}'}\right)$ and $\bar{k}' = \sqrt{1-\bar{k}^2}$, and $1 > \bar{k} > 0$ and $\dfrac{\pi}{2} > \psi > 0$.

### 4.8.3  Integrals Containing $\ln\left(\dfrac{\varepsilon+u}{\varepsilon-u}\right)$ and Series Representation

(4.51) $\int_\alpha^\beta \ln\left(\dfrac{\varepsilon+u}{\varepsilon-u}\right)\dfrac{u^2 \, du}{\sqrt{(u^2-\alpha^2)(\beta^2-u^2)}} = \dfrac{\pi}{\varepsilon}\left(\varepsilon^2 - \sqrt{(\varepsilon^2-\alpha^2)(\varepsilon^2-\beta^2)}\right) + \pi\beta\left[F(\varphi,k) - E(\varphi,k)\right]$,

where $\varphi = \arcsin\dfrac{\beta}{\varepsilon}$ and $k = \dfrac{\alpha}{\beta}$, and $\varepsilon > \beta > \alpha > 0$.

(4.52) $\int_{e_2}^{e_1} \ln\left(\frac{1+q}{1-q}\right) \frac{dq}{\sqrt{(e_1^2 - q^2)(q^2 - e_2^2)}} = \pi \sum_{m=0}^{\infty} \Lambda_{2m+1}$,

(4.53) $\Lambda_1 = 1$,

(4.54) $\Lambda_3 = \frac{e_1^2 + e_2^2}{6}$,

(4.55) $\Lambda_{2n+5} = \frac{(e_1^2 + e_2^2)(2n+3)^2 \Lambda_{2n+3} - 2e_1^2 e_2^2 (n+1)(2n+1)\Lambda_{2n+1}}{2(n+2)(2n+5)}$,

where $n = 0, 1, 2, 3, \ldots$, $m = n + 2$, and $1 > e_1 > e_2 > 0$.

(4.56) $\int_{e_2}^{e_1} \ln\left(\frac{1+q}{1-q}\right) \frac{q^2 \, dq}{\sqrt{(e_1^2 - q^2)(q^2 - e_2^2)}} = \pi \left( \sum_{m=0}^{\infty} \Omega_{2m+1} - 1 \right)$,

(4.57) $\Omega_1 = 1$,

(4.58) $\Omega_3 = \frac{e_1^2 + e_2^2}{2}$,

(4.59) $\Omega_{2n+5} = \frac{(e_1^2 + e_2^2)(2n+1)(2n+3)\Omega_{2n+3} - 2e_1^2 e_2^2 (n+1)|2n-1|\Omega_{2n+1}}{2(n+2)(2n+3)}$,

where $n = 0, 1, 2, 3, \ldots$, $m = n + 2$, and $1 > e_1 > e_2 > 0$.

## 5    Concluding Remarks

Motivated by a study of the surface area of arbitrary ellipsoids, a number of potentially useful integrals were derived and evaluated. The integrals $I_N(a,b,c)$, $N = 1, 2, \ldots, 6$, and their analytic extensions, naturally extend the integral tables found in the literature. Of the twelve integrals and extensions tabulated and described in Sections 4.7 and 4.8, only three were directly related to integrals found in the literature, while two other extensions could be converted to integrals found in the literature. Although the integrals arise in a straightforward derivation of the surface area of arbitrary ellipsoids, they are expected to be useful in other applications.

Several techniques were employed to verify the correctness of the integrals $I_N(a,b,c)$ and their analytic extensions. References were provided, wherever possible, throughout the derivation and evaluation of integrals and extensions; analytic-extension interrelationships were studied, some analytic extensions were compared with standard references, and the numerical equivalence of the integrals and extensions was studied by the repeated application of Mathematica 5 within the stated parameter space. Related infinite series representations of two ellipsoid-specific integrals, namely $\Sigma_1(a,b,c)$ and $\Sigma_2(a,b,c)$, were similarly derived and evaluated, and the correctness of $\Sigma_1(a,b,c)$ and $\Sigma_2(a,b,c)$, similarly verified.

Finally, the entire investigation provides an exercise in integration techniques involving elliptic integrals and the surface area of arbitrary ellipsoids that could prove useful and stimulating to students and professionals alike.

# 6 Acknowledgment

The authors acknowledge the substantial contributions made to this study by Mathematica 5. In particular, the numerical equivalence of the nine integrals and two series presented in Section 4.8 were studied and tested by an extensive application of Mathematica 5 within the stated parameter space. In addition, the graphics in Figs. 1 and 3 were created using Mathematica 5. The authors gratefully acknowledge assistance with the general references provided by Ronaele Freestone and Susan Heckethorn at the Los Alamos National Laboratory Research Library.

# 7 General and Bibliographic References

Two separate but related reference lists are provided. General References provide additional information often of a historical nature that may be useful to the reader. Bibliographic References were used in the production of this document. Their application is denoted by bold type within the text of the document and at the right margin of equations on the appropriate page.

## 7.1 General References

## 7.2   Bibliographic References